\documentclass{amsart}
  \usepackage{amscd,amssymb,epsfig}
   \usepackage{epic,eepic}

                     \numberwithin{equation}{subsection}

                     \newtheorem{propo}{Proposition}[subsection]
                     \newtheorem{corol}[propo]{Corollary}
                     \newtheorem{theor}[propo]{Theorem}
                     \newtheorem{lemma}[propo]{Lemma}
                     \theoremstyle{definition}
                     
                     \newtheorem{examp}[propo]{Example}
                     \theoremstyle{remark}
                     \newtheorem{remar}[propo]{Remark}

		     \newcommand{\CC}{\mathbb{C}}
		     \newcommand{\QQ}{\mathbb{Q}}
                     \newcommand{\ZZ}{\mathbb{Z}}
                     \newcommand{\RR}{\mathbb{R}}

		     \newcommand{\Ker}{\operatorname{Ker}}

\newcommand{\Fl}{\operatorname{Fl}}
\newcommand{\cor}{\operatorname{cr_1}}
\newcommand{\cro}{\operatorname{cr}}
\newcommand{\Iso}{\operatorname{Iso}}

                     \newcommand{\id}{\operatorname{id}}

		    \newcommand{\modu}{\operatorname{mod}}
\newcommand{\Int}{\operatorname{Int}}
		     
		  \newcommand{\rank}{\operatorname{rank}}

                    \newcommand{\card}{\operatorname{card}} 
 \newcommand{\cross}{\operatorname{cr_{>1}}}
 \newcommand{\Cross}{\operatorname{sing}}
\newcommand{\Aut} {\operatorname {Aut}}
\newcommand{\Auut} {\operatorname {Aut}}

\newcommand{\Stab} {\operatorname {Stab}}
 
\begin{document}
      \title{Curves on surfaces, charts, and words}
                     \author[Vladimir Turaev]{Vladimir Turaev}
                     \address{%
              IRMA, Universit\'e Louis  Pasteur - C.N.R.S., \newline
\indent  7 rue Ren\'e Descartes \newline
                     \indent F-67084 Strasbourg \newline
                     \indent France }
                     \begin{abstract} 	 We give a combinatorial description of closed curves on   oriented surfaces in terms of
certain permutations, called charts. We   describe automorphisms of curves in  terms of charts and compute the total number of curves
counted with appropriate weights. We also discuss   relations between curves, words, and complex structures on surfaces.
                     \end{abstract}
                     \maketitle

                  \section{Introduction} 
	 
Words, defined as finite sequences of letters,   were   used to describe geometric figures   already by the Greeks. It was common to
label points by  letters and to encode   polygons by   sequences of labels of their vertices, see \cite{eu}.  Gauss \cite{ga} extended this
method  to    closed curves with self-intersections.  He considered   generic  curves     on the plane that  is  smooth  immersions $S^1\to
\RR^2$  with    only double transversal self-intersections.  The  Gauss   word  of  a generic curve  is  obtained by  
labelling its  self-crossings  by different letters and writing down  these letters in the  order of  their appearance when
one traverses   the  curve.  Each letter appears in this word twice.  Since the  work of Gauss, considerable efforts  were directed at
characterizing the words   arising  in this  way from generic curves on $\RR^2$ and on other surfaces, see
\cite{ro},
\cite{rr},
\cite{dt},
\cite{cw}, \cite{ce}, \cite{cr}.

The present paper  is motivated by the   desire to invert the procedure and to find a geometric presentation  of arbitrary  words   by  
curves  on  (oriented) surfaces.  A related question is to find   a   topological classification  of   curves on  oriented surfaces.   We consider 
in this paper  all (closed) curves   with finite number of self-crossings of any finite multiplicity
$\geq 2$.   We endow curves  with a finite set of distinguished
points of multiplicity 1,  called   corners.   The set of corners  is choosen arbitrary and  may be empty.   The  Gauss    word
extends to such curves but   is far from being  a full   invariant.

   We introduce an     invariant
of    curves   called the \lq\lq chart". The chart of a curve  $f$  is   a  permutation of the set
$ \{-n, -n+1,..., 
-1,1, ..., n-1, n\}$ where  $n\geq 1$ is an integer determined by  $f$.
  The definition of the
chart  of a curve reminds of  Grothendieck's cartographic groups, cf. 
\cite{js}.  

 Knowing the chart of a curve,   we can   recover  the position of the curve  in its  narrow neighborhood in the ambient surface.   For  a
{\it filling} curve, that is a curve   on a closed oriented surface  whose complementary regions     are disks,  the
chart determines  the curve  and the surface up to homeomorphism.  In particular,    the genus of the surface can be computed directly
from the chart;  we give a combinatorial formula for this genus.

A curve $f:S^1\to \Sigma$ may have  automorphisms  that is   degree 1  homeomorphisms $\Sigma\to \Sigma$ preserving
both the curve  and
 its set of   corners.  The   isotopy
 classes of  such automorphisms    form a group  $\Auut(f)$.  If $f$ is    filling, then $\Auut(f)$
  is a finite cyclic group.  We compute its order    
$\vert
\Auut (f)\vert$  in terms of the  chart  of $f$. 

 For self-transversal curves, i.e., curves  with only transversal self-intersections (of arbitrary multiplicity),
 the information encoded in the chart can be presented in a
more compact   form of a   semichart. A semichart   is a pair (a permutation of the set $\{1,2,...,n\}$, a subset of
this set). Filling self-transversal curves  are determined by their semicharts  up to homeomorphism.

We use  charts to   compute the number of     curves counted 
 with certain weights. We state here a version of this result      for  filling  self-transversal  
curves  (on closed oriented surfaces  of arbitrary genus).  Let  $\mathcal C_{str}$ be  the set of homeomorphism classes of such
curves.  Then we have an equality of formal power series
\begin{equation}\label{formal2}\sum_{f\in \mathcal C_{str}}
\,\frac { 1 } {(n(f)-1) !\,
\vert
\Auut(f)\vert} \,\, t_1^{k_1(f)} t_2^{k_2(f)} t_3^{k_3(f)}\cdots   =     \prod_{m\geq
1}
\exp 
 \left ( \frac{2^{m-1}}{m} t_m   \right )     \end{equation}
where    $t_1,t_2,...$ are independent
variables,  $k_1(f)$ is the number of corners of $f$,  $k_m(f)$ with $m\geq 2$ is the number of crossings of $f$ of multiplicity $m$,
and $n(f)=\sum_{m\geq 1} m\, k_m(f)$. The  monomials on the left hand side of (\ref{formal2}) are finite since $k_m(f)=0$ for  any $f$
and all sufficiently big
$m$. Formula \ref{formal2} indicates that  in  appropriate statistical terms, the crossings  
 of different
multiplicities are independent of each other.

Comparing the coefficients of the monomial $t_1^{k'} t_2^k$  with $k',k\geq
0$ on both sides  of (\ref{formal2}), 
we obtain  
\begin{equation}\label{eefr25} \sum_{f\in \mathcal C_{str} (k',k)} \frac {1} {\vert \Auut  (f)\vert}  
=  \frac{( 2k+k' -1)!}{k!\, k'! }  \end {equation}
where $\mathcal C_{str} (k',k)\subset  \mathcal C_{str} $ is the set   of homeomorphism  classes of filling     curves with $k'$
corners,
$k$  double transversal  crossings, and no other crossings. In particular,  for $k'=0$,  
$$ \sum_{f\in \mathcal C_{str}(0,k)} \frac {1} {\vert \Auut  (f)\vert}  
=  \frac{(2k -1)!}{ k! }  .$$
For $k'=1$, Formula \ref{eefr25} gives
  $\card (\mathcal C_{str} (1,k))
=  2k!/k!$  since  $ \Auut  (f)=1$ for   $f\in  \mathcal C_{str} (1,k)$.  

 Although   curves are   purely topological
objects, they are related  to deep algebra and geometry.   The image of a filling curve $f:S^1\to \Sigma$ is a  graph whose vertices are 
the self-intersections and the corners of $f$. In terminology of
\cite{sc}, p. 51,  this graph is a pre-clean dessin d'enfants.  By the classical Grothendieck construction, it 
induces on $\Sigma$ a structure of an algebraic curve over
$\overline {\QQ}$.  Over  $  {\CC}$ this gives a smooth  curve  with distinguished set of  $k=\sum_{m\geq 1} k_m (f)$ points 
consisting of the self-intersections and the corners of $f$. Numerating them, we obtain a   point of the moduli space $\mathcal M_{g,k}$
where $g=g(\Sigma)$ is the genus of
$\Sigma$. 

Another  geometric construction applies when $\Sigma$ is  a smooth   surface.   It associates with a curve $f:S^1\to \Sigma$ an
oriented  
  knot  in  the 3-manifold $S\Sigma$ formed by
  tangent vectors of $\Sigma$ of length 1 with respect to a certain Riemannian metric on $\Sigma$.  If   $f$
is smooth and self-transversal,   this knot    is formed by the unit 
positive tangent   vectors of $f$.    This knot is transversal to the standard contact structure on $S\Sigma$. For general $f$,   
  we   first approximate it by a smooth self-transversal curve and then proceed as above.

These geometric constructions   can be combined with the construction of filling curves from
charts and allow to associate with each chart   an algebraic curve over
$\overline {\QQ}$,  a point of a moduli space, and   a knot in the tangent circle bundle over a surface. 
 
Observe also that the charts of self-transversal curves can be naturally viewed as elements of   Coxeter groups  of type $B$.

Coming back to words,  note that  for  curves on oriented surfaces, the letters in the associated  word naturally acquire signs $\pm$.
 To simplify notation, we    omit the minuses.  Thus instead of writing   $ A^-B^+A^+C^-$ we write
$ AB^+A^+C$. We show that every such word (in any  alphabet)  can be realized by a curve on a surface. This
realization is by no means unique.  Not all words can be realized by  self-transversal curves. We say
that a word   is odd if every  letter appearing in this word (with or
without   superscript $+$)   appears   without    + an odd number of times.  For instance,  the words  $  ABC, AA^+, A^+B^+AB$ are  odd
while  
  $  AA, A^+A^+, A B^+$  are not odd. We prove that    a  word
  can be realized by a self-transversal curve if and only if  it is odd. 

A deeper connection  between words and curves involves so-called coherent curves, generalizing the generic curves. Note that the
branches  of  a self-transversal curve $f:S^1\to \Sigma$ passing through a given crossing     acquire a cyclic order  obtained by moving on
$\Sigma$ around  the crossing.  We call $f$   coherent if  for all its crossings, 
  this  cyclic order  coincides with the   order of appearance  of these  branches when one traverses the curve. 
We prove that  any odd word has a unique realization  by a  coherent   filling (self-transversal) curve.  This gives a   geometric
form to odd words: they correspond bijectively  to  the  homeomorphism classes of    coherent filling curves.  For example, the   word  
$ABC$   corresponds to a  triangle on $  S^2$ with corners in the vertices. More generally,  the word 
$A_1A_2...A_n$  formed by  
  distinct letters
$A_1, A_2, ...,A_n$  corresponds to an embedded
$n$-gon 
 on $ S^2$ with corners in the vertices.  The   word $AA^+$ corresponds to a 8-like curve on $S^2$. The word $A^+B^+AB$
corresponds to a  curve on the torus.  We formulate  conditions on an odd word necessary and sufficient for the
corresponding coherent curve to be planar.  This extends  a theorem of Rosenstiehl \cite{ro} on Gauss words.

Combining with the geometric constructions outlined above, we   associate with each odd word   an algebraic curve over
$\overline {\QQ}$,  a point of a moduli space, and   a knot in a 3-manifold.

We also  introduce and study  further geometric classes of curves  (pointed/alternating/beaming/perfect curves).    All
these classes can  be described in terms of their charts.

 Throughout the paper, by   a
surface, we mean an {\it oriented} 2-dimensional manifold possibly with boundary.  

\section{Curves and homeomorphisms}

  \subsection{Tame curves}\label{loca} 
 To avoid locally wild behavior, we  consider only    tame curves. 
A {\it tame curve} on a surface $\Sigma$ 
is  a continuous map  $f$ from   $ S^1= \{z\in \CC\,\vert\, \vert z\vert = 1 \}$  to  $\Sigma - \partial
\Sigma$ satisfying   two conditions:

(i) $f$ is locally injective, i.e., each point of $S^1$ has a
  neighborhood $U$ such that
$f\vert_U:U\to \Sigma$ is injective;

(ii)  for any      $a\in f(S^1)$ there are a closed neighborhood $V\subset \Sigma$ and a
homeomorphism of $ V $ onto the unit complex disk $D= \{z\in \CC\,\vert\, \vert z\vert \leq 1 \}$ sending $a$  to $0$
and sending 
$   f(S^1)\cap V$ onto the set $\{z\in D\,\vert\, z^m\in \RR \}$ where  $m $ is a positive integer (depending on $a$).

The  integer $m=m_a $   in (ii) is  the {\it
multiplicity} of $a$. The    set $\{z\in D\,\vert\, z^m\in \RR \}$  consists of
$2m$ radii of  $D$. Moving along a radius  towards $a=0$ the curve  then goes away  along another radius.  
Clearly,   $\card (f^{-1}(a))=m$.   

A  
{\it crossing}   of a  tame curve
$f:S^1\to
\Sigma$    is a point of $ f(S^1)$ of multiplicity
$>1$.  
The set of crossings of  $f$ is
denoted $\cross (f)$.    The compacteness
of
$S^1$ implies  that    this set   is  finite. Observe that
$$f^{-1} (\cross (f)  )=\{x\in S^1\,\vert \,\, {\rm {there\,\,
is}}\,\,  y\in S^1- \{x\} 
\,\,  {\rm {such
\,\, that }}\,\,  f(x)=f(y)\}.$$   

 If $\Sigma$ is a smooth surface, then    all  smooth immersions $S^1\to
\Sigma$ with finite number of crossings are tame.

\subsection{Corners and filling curves}\label{regu}   A  {\it tame 
  curve with corners} on a surface
$\Sigma$ is a tame curve $f:S^1\to
\Sigma$ endowed with a finite  (possibly empty) subset  of  $ f(S^1)-  \cross (f)$. This finite subset is denoted  $\cor (f)$, its points  are
called   {\it corners} of  $f$.  They all have multiplicity 1.   
Set   $\cro (f)= \cor (f)\cup \cross (f)$ and $ \Cross(f)=f^{-1} (\cro (f))$.
Both these sets are finite.  
Set   
$$n(f) =\card (\Cross (f))=\sum_{a\in \cro (f)} m_a  =\sum_{m\geq 1} mk_m (f)$$
where  
$k_1(f)$ is the number of corners of
$f$ and    $k_m(f)$ with $m\geq 2$  is the number of crossings  of
$f$ of multiplicity $m $.  Clearly,   $  n(f) = 0 $  if and only if $f$ has no crossings
and no corners.

\vskip0.2truecm
{\it From now on, the word \lq\lq curve"  means a  tame  curve with corners on a  surface.}
\vskip0.2truecm
      A 
 curve $f:S^1\to \Sigma$ is {\it filling} if $\Sigma$ is  a closed connected   surface 
and 
  all components of
$\Sigma-  f(S^1)$ are  open 2-disks.   
Any curve $f$ on any surface $ \Sigma$   gives rise to a filling  
curve by taking a  regular closed   neighborhood $U\subset \Sigma$ of $f(S^1)$ 
and gluing   2-disks to the components of $\partial U$.

 \subsection{Homeomorphisms}\label{homeo}  A {\it homeomorphism}   of  curves $f_1:S^1\to
\Sigma_1,  f_2:S^1\to \Sigma_2$ is a pair  $(\varphi^{(1)}: S^1 \to S^1, \varphi^{(2)} :\Sigma_1 \to \Sigma_2)$
of orientation  preserving  
  homeomorphisms   such that
$f_2\, \varphi^{(1)}=\varphi^{(2)} f_1  $ and
$\varphi^{(2)}(\cor (f_1))=\cor (f_2)$.
If there is such  a pair, then  $f_1 $
and $ f_2 $ are {\it homeomorphic}.  The  role of $\varphi^{(1)}$  is to ensure that    curves obtained from each other 
by  re-parametrization are homeomorphic.  Note that $\varphi^{(1)} $ maps 
$\Cross (f_1)$  bijectively onto
$\Cross (f_2)$ and  $\varphi^{(2)}$ maps    $f_1(S^1), \cross (f_1), 
\cor (f_1)$  bijectively onto   $f_2(S^1), \cross (f_2),   \cor (f_2)$, respectively.  

  A curve homeomorphic to a
filling curve is itself filling.  The main problem of the combinatorial  theory of curves   is a classification of  filling 
  curves up to homeomorphism.

\subsection{Examples}\label{exam} 1.    An embedding
$S^1\hookrightarrow S^2$ with empty set of corners is a {\it trivial} curve. Any two trivial curves are homeomorphic.

2.  Let  
$A,B,C,D\in \RR^2\subset S^2$ be  consecutive  vertices of a square with center $O$
  where $\RR^2$  is oriented so that the pair of vectors $(AB,AC)$ is positive.  The closed broken lines  $ABODCOA$  and 
$ABOCDOA$ with  no  corners   are  filling  curves on $S^2$ with 
 one crossing   $O$ of multiplicity  2.  They are not homeomorphic.

 3.  Let  
$A,B,C,D,E,F\in \RR^2\subset S^2$ be  consecutive  vertices of   a regular   hexagon   with center $O$
 where the   plane is oriented as above.   The  closed broken line
$ABOEFOCDOA$  with no corners  is a filling curve on $S^2$ (called trifolium)  with one crossing $O$ of multiplicity 3.

\section{Charts and flags}
                    
 \subsection{Preliminaries}\label{prels} For a positive integer $n$, set
 $\hat n=\{1,2,...,n\}$ 
and
$$\overline n=(-\hat n) \cup \hat n= \{-n, -n+1,...,-2,
-1,1,2,..., n-1, n\}.$$
The {\it circular permutation}
$\sigma_n:\overline n \to \overline n$ sends $ \pm k$ to $\pm (k+1)$ for $k=1,2,..., n-1$ and sends
$ \pm n$ to $ \pm 1$. For $n=0$,   set  $\hat n =\emptyset$, $\overline n=\{0\}$, and $\sigma_0=\id:\{0\}\to \{0\}$.

For a  finite set  $F$  and a bijection $t:F\to F$, an {\it orbit} of  $t$, or shorter a {\it $t$-orbit},  is a minimal non-empty
$t$-invariant subset of
$F$. 
  The unique orbit of $t$ containing a given element $a\in  F$ consists of $m$ elements $ a, t(a), t^2(a), ..., t^{m-1}(a) $ where 
$m $ is the minimal positive integer such that $t^m(a)=a$.  The mapping $t$ determines a      cyclic order $\prec$ 
on this  orbit by  $ a\prec t(a) \prec  t^2(a)\prec ...\prec t^{m-1}(a)\prec  a $. It is clear that $F$ is a disjoint
union of   orbits of
$t$. The set of    orbits of
$t$ is denoted 
$F /t$.

  \subsection{Charts}\label{charts} A   pair $(n,t)$ consisting of
an integer 
$n\geq 0$   and  a bijection
$t:\overline n \to \overline n$ is  a  {\it chart} if      for any
$k\in
\overline n$, there is   $s\in \ZZ$  such that
$ t^s(k)=-k$.  The  chart  $(0,\id:\{0\}\to \{0\})$  is called the {\it trivial} chart. 
In a non-trivial chart $(n,t)$   the bijection  $t$ is necessarily
fixed-point free: if 
$t(k)=k$, then  
$t^s(k)=k\neq -k$ for   $s\in \ZZ$. Any orbit of $t$ is  invariant under the negation $\overline n \to \overline n, k\mapsto -k$
and   has an even number of elements.
The cyclic group $\ZZ/n\ZZ$ acts on the set of charts $(n,t)$ via     $t\mapsto \sigma_n t (\sigma_n)^{-1}$.  

For each curve $f$, we shall define a chart $(n(f),t(f))$ where $ n(f)$ is the number defined in Section \ref{regu}. The following theorem
gives a  topological classification of filling  curves    in terms of charts. 

\begin{theor}\label{t5} The formula $f\mapsto (n(f),t(f))$  defines a bijective correspondence between    filling   curves  considered up to
homeomorphism  and   charts  
  considered up to conjugation by the   circular permutation.
   \end{theor}

Theorem \ref{t5} directly follows from Lemmas \ref{t2}  and \ref{t4} stated below.  This theorem   provides  a   combinatorial 
method of presenting filling  curves: to specify a curve, it suffices to specify its chart.

Note that all  curves on   $S^2$ are filling. Therefore the chart   is a  full topological invariant of   curves
on
$S^2$.

      \subsection{Flags}\label{flags}    
 A {\it flag}  
of  a curve  $f $  is a pair   ($x\in  \Cross (f)  , \varepsilon=\pm  $).   The flag $(x,-)$ is     {\it
incoming},  the flag $(x,+)$ is   {\it outgoing}.   The flags
$(x,-)$ and 
$(x,+)$ are  said to be  {\it opposite}.  
With a flag
$r= (x,\varepsilon) $ we associate a small arc in
$S^1$  beginning at
$x$ and going counterclockwise if $\varepsilon=+$ and clockwise if $\varepsilon=-$. 
The image of this arc under $f$ is a small embedded arc on
$f(S^1)$ with one endpoint  at
$f(x)$.  The latter arc presents $r$ geometrically.  We say that $r$ is a {\it flag    at   }  $f(x)$  and   $f(x)$    is  
the {\it root} of
$r$.  Each
$a\in
\cro (f)$ is a root of
$2m_a$ flags. 

The set of flags of $f$ is denoted  $\Fl
(f)$.  This set     has $2 n$ elements  where  $ n =n(f)=\card (\Cross (f))$. If $n\neq 0$, then $\Fl (f)$   can be identified
with the set   $\overline n $ as follows. Starting at a point of
$S^1- \Cross (f)$ and  traversing $S^1$  counterclockwise  denote the points of
$\Cross (f)$ consecutively   $x_1, x_2,...,x_n$. The  flag $(x_k,\pm)$   corresponds to  $\pm\, k \in \overline
n$ for  $k=1,2,...,n$. The resulting  bijection $\Fl(f)\approx  \overline n$   is well defined up to composing  with   
$\sigma_n$.

 \subsection{Charts of curves}\label{chartscurv}   For a curve $f:S^1\to
\Sigma$, we   define  a  bijective map   $t: \Fl (f)\to \Fl (f)$ called {\it   flag rotation}.    Starting from a   flag of $f$ at  $a\in \cro (f)$
  and circularly moving on $\Sigma$  in the positive direction around  
   $a $, we   numerate   the flags    of $f$
at $a$ 
  by the residues $1,2,..., 2 m_{a}-1, 2 m_{a}\, (\modu 2m_a)$. The map  $t$  transforms the $k$-th flag at 
$a$  into the $(k+1)$-st flag   at  $a$.     In particular, if $a\in \cor(f)$, then $t$
permutes the two (opposite) flags at 
$a$.

   Assigning to each flag its root, we obtain  
 $\Fl (f)  /t= \cro (f)$.
  The  orbit of $t$  
 corresponding to $a\in \cro (f)$   consists of  the $2m_a$ flags of $f $ at  $a$.  Opposite flags   always  lie  in the same  orbit.
 
If  $n (f)=0$, then by definition the chart of $f$ is  the trivial chart  $(0, \id :\overline 0\to
\overline 0)$.  Suppose that    $n=n (f)\neq 0$. Conjugating $t: \Fl (f)\to \Fl (f)$ by  the  bijection  $\Fl(f)=\overline n$, we obtain  a
bijection
$\overline n
\to \overline n$ also denoted $t$ or $t(f)$. It is
determined  by $f$   up to   conjugation by $\sigma_n$. The pair $(n,t:\overline n
\to \overline n)$  is a chart called 
the {\it chart of}  $f$.  Note the  identifications $\overline n/t=\Fl(f)/t=\cro(f)$.   The opposite choice of orientation on
$\Sigma$     yields  
$(n, t^{-1})$.    Since the   chart of    $f$  is entirely  defined   inside  a narrow neighborhood of $f(S^1)\subset
\Sigma$,  it depends only on the  filling curve determined by $f$.

 	\begin{lemma}\label{t2} Any  chart  $(n,t)$  is  the chart of a filling  curve.
   \end{lemma} 

  \begin{proof}     It is enough to consider the case
$n\geq 1$. Set  
$x_k=\exp (2k
\pi i /n)\in S^1$ for  $k\in \hat n= \{1, 2, ...,n\}$. Starting in $x_k$ and moving   along $S^1$
counterclockwise (resp. clockwise) to the distance
$\pi/2n$ we sweep  an arc  on $S^1$  denoted 
$\alpha_{k} $ (resp. $\alpha_{-k} $).   The  arcs $\{\alpha_{k}  , \alpha_{-k} \}_{k\in \hat n}$ are disjoint except that  $\alpha_{k}\cap
\alpha_{-k}=\{x_{k}\}$ for all $k$.

  Identifying $x_k$ and $x_j$ (with $k,j\in \hat n $) each time
there is a power of $t$   transforming  $k$ into $j$, we obtain from 
$S^1$   a 1-dimensional CW-complex  
$\Gamma$.    The 0-cells (vertices) of $\Gamma$ are the 
points 
$p(x_1), p( x_2),..., p(x_n)$ where $p:S^1\to \Gamma$ is  the natural projection.  All    other points of $\Gamma$ form the 
open  1-cells.    

 We  shall  thicken
$\Gamma$ to a surface. First we thicken each
vertex 
$v\in \Gamma$   to a disk as follows.  Set   $J_+=\{j\in \hat n,\vert\, p(x_j)=v\}$ and   
$J=(-J_+)\cup  J_+\subset \overline n$.  We claim that $J$ is an orbit of  $t$. 
 To see this, pick any $k\in J_+$ and denote its $t$-orbit by $[k]$. We show that $J=[k]$. By the
definition of
$\Gamma$, we have $J_+\subset [k]$.  Hence, by the definition of a chart,      $-J_+ \subset [k]$ and  $J\subset [k]$.
To prove the  inclusion $[k]\subset J$, we show that 
$t^s k\in J$ for  any  $ s\in \ZZ$.  Pick $s_-\in \ZZ$ such that
$-t^s k=t^{s_-}  k  $.  If
$t^s k\in \hat n$, then
$t^s k\in J_+ \subset  J
 $ by the definition of $\Gamma$.  If    $t^s k\in -\hat n$, then   $-t^s k=t^{s_-}  k \in J_+$ and $t^s k\in -J_+ \subset  J$.

 Set $m=\card (J_+)$.
The set $V_v=
\cup_{r\in J}
\alpha_r$ is a neighborhood of
$v$ in
$\Gamma$ consisting of
$2m$  embedded arcs  $\{  \alpha_{r} \}_{r\in J} $
meeting   at $v$.  We
embed
$V_v$ in a copy $D (v)$ of  the unit 2-disk
$D=
\{z\in
\CC\,\vert\,
\vert z\vert
\leq 1
\}$ as follows. Pick   $r\in J$. The embedding $V_v\to  D(v)$
sends
$v$  to
$0$ and sends the arc $\alpha_{t^s r}$ onto the radius $R_s= \exp ( s \pi i  /m) \cdot [0,1]$ of $D(v)$ for  $s=0,1,..., 2m-1$. 
 In this way $v$ is thickened 
 to a copy $D(v)$ of $D$.   We  endow $D(v)$  with counterclockwise orientation.   
Note that going around $v$ in the positive direction we cross  the arcs $\{  \alpha_{r} \}_{r\in J}$  in the cyclic order
determined by
$t$  as in Section \ref{prels}.

Endow a  1-cell  $\gamma$ of $\Gamma$   with the orientation induced by the counterclockwise orientation on
$S^1$.  Then $\gamma$ leads from a vertex $v_1$ of $\Gamma$ to a vertex $v_2$ of  $\Gamma$ (possibly
$v_1=v_2$). These  vertices 
are  thickened  above to disks
$D_1=D(v_1), D_2=D(v_2)$.  For $l=1, 2$, the disk $D_l$    meets $\gamma$     along a
radius 
$R_{s(l) } $  as above with endpoint $a_l \in  \partial D_l =S^1$.  The part of
$\gamma$ not lying in $\Int D_1\cup \Int D_2$ is a closed  interval $\gamma_\bullet\subset  \gamma$ with endpoints  
$a_1, a_2$. We thicken  
$\gamma_\bullet $ to a ribbon 
$\gamma_\bullet \times  [- 1, 1]$  glued to  
$D_1\cup  D_2$  along  the    embedding   $a_1\times  [- 1, 1] \hookrightarrow \partial D_1$  
sending   $(a_l, u\in  [- 1, 1] )$ to $\exp (u\pi i /2n)  \cdot a_1 $
and the  embedding   $a_2\times  [- 1, 1] \hookrightarrow \partial D_2$  
sending   $(a_2, u )$ to $\exp (- u\pi i /2n)  \cdot a_2 $. Note that     the   orientation   on   $D_1, D_2$ extends to
 their union with the ribbon.

Thickening in this way  all   cells of $\Gamma$ we embed $\Gamma$ into a 
surface  $ 
U$.  By   
construction,  $U $ is a   compact connected  oriented   
surface with non-void 
boundary. All components
 of
$U -  \Gamma$ are  homeomorphic to   $S^1\times [0,1)$.    Gluing  2-disks to the components of $\partial
U$ we obtain   a closed connected  oriented surface
$\Sigma\supset U\supset \Gamma$ such that all components of $\Sigma -  \Gamma$ are open disks.    Composing the  
projection
$p: S^1\to
\Gamma$ with the  inclusion $\Gamma \hookrightarrow \Sigma $ we obtain a filling   curve   $ S^1\to \Sigma $. 
We choose its set of corners to be  $\{ p(x_k)\,\vert \, k\in \hat n,  t(k)=-k\}$. It is
clear from the construction that the chart of the resulting  curve is the original   chart $(n,t)$. 
 \end{proof} 

 \begin{remar}\label{figjgig}  The
construction of   $\Sigma $ is well known in the context of generic curves, see 
\cite{f}, 
\cite{ca}, \cite{cw}.  \end{remar}

 \subsection{Action of homeomorphisms}\label{action}  A homeomorphism $\varphi=(\varphi^{(1)}: S^1 \to
S^1,
\varphi^{(2)} :\Sigma_1 \to \Sigma_2)$ of   curves $f_1:S^1\to
\Sigma_1,  f_2:S^1\to \Sigma_2$  induces a bijection
$\varphi_*:\Fl(f_1)
\to
\Fl(f_2)$.  It sends the flag of
$f_1$ represented by   a small embedded arc  
$\alpha\subset f_1(S^1)$ with   root   
$a\in \cro (f_1)$   to the flag of $f_2$ represented by the arc $\varphi^{(2)} (\alpha)\subset f_2(S^1)$
with   root   
$\varphi^{(2)}(a)\in \cro (f_2) $.     The same   $\varphi_*$  is defined by 
$\varphi_*((x, \pm))= (\varphi^{(1)} (x), \pm)$ where   $x\in  \Cross (f_1)$.  The bijection  $\varphi_*$ 
  sends opposite flags to opposite flags and commutes with the flag rotation.  
This implies   that  homeomorphic  curves   have the same
charts.  
The next
lemma establishes the converse for    filling   curves.

	\begin{lemma}\label{t4}   Filling   curves having the same  charts are
homeomorphic.
   \end{lemma} 
\begin{proof}   If one of the curves is trivial, then the claim is obvious.  Consider the case of non-trivial curves.  The
reconstruction of a curve from its chart given in the proof of  Lemma
\ref{t2} is   canonical except at the place where  we pick $r \in J$.   A
different choice of
$r$  would lead  to another embedding
$V_v\hookrightarrow  D$  obtained  from the first one   by composing with a   rotation of $D$ around its center  to an
angle proportional to $2\pi/m$. This, however,  gives  the same thickening of $v$ to a  disk $D(v)$; only the identification 
of  $D(v)$ with
$D$   differs by this rotation.    Thus,  
 knowing the chart of a  curve $f:S^1\to \Sigma$ we can  
reconstruct a regular neighborhood  $U\subset \Sigma$   of $f(S^1)$   and the curve $f:S^1\to U$  up to
homeomorphism.  Since a homeomorphism of circles extends to a homeomorphism of disks bounded by these circles (cf.
Lemma \ref{lemop1} (i)  below), we conclude   that knowing the chart of a  filling   curve   we can  
reconstruct   the curve  up to
homeomorphism.   \end{proof}

 \subsection{Examples}\label{exa} 1. The chart of an embedding  $ S^1\hookrightarrow  S^2$    with $n $ corners  is the
negation    $\overline n \to \overline n, k\mapsto -k$.

2. The chart of  the  curve   $ABOCDOA$ (Example \ref{exam}.2)   is the cyclic permutation 
$(1,-2,2,-1)$  of $\overline 2$  sending $1$ to $-2$, $-2$ to $2$, $2$ to $-1$, and
$-1$ to
$1$.  The chart of  the  8-like curve  $ABODCOA$  (the same example)  is the cyclic permutation 
$(1,2,-1,-2)$  of $\overline 2$.   The chart of the same  curve  $ABODCOA$  with corners $A,B,C,D$ is   the permutation $(1,-1)
(2,-2) (3,6,-3,-6) (4,-4) (5,-5)$ of $\overline 6$.

 3.  The chart of  the  curve 
$ABOEFOCDOA$  (Example \ref{exam}.3) is the cyclic permutation 
$(1,-2, 3, -1,2,-3)$  of $\overline 3$.

  \subsection{Riemann surfaces}\label{encoreeee}   By the results above,  any  chart
$(n,t)$ determines  a filling curve $f:S^1\to \Sigma$ uniquely up to homeomorphism. Applying to the graph $f(S^1)$
the Grothendieck construction,  we obtain  an algebraic curve over $\overline {\QQ}$
and   a  point of the moduli space
$\mathcal M_{g,k}$ where $g=g(\Sigma)$ and $ k=\card (\cro (f))= \card (\overline n/t)$. Here to fix an order on the set
$\cro (f)$ we identify it   with the set of orbits of  
$t:\overline n
\to
\overline n$ and order  the latter      by 
$A<B$ if 
$\min_{a\in A}
\vert  a\vert < \min_{b\in B} \vert  b\vert$.    
  Note that under conjugation of
$t$ by $\sigma_n$,  this order   may change.  

These constructions can be generalized in various directions. Having  positive real numbers $r_1,...,r_n$ we can provide the graph
$f(S^1)$ with a metric   assigning $r_i$ to the only edge containing  the
$i$-th outgoing flag  for $i=1,...,n$.  The resulting  metric graph determines a complex structure on $\Sigma$, see for instance \cite{mp}.
For
$r_1=...=r_n=1$, this gives   the  same Riemann surface as in the Grothendieck  construction on $f(S^1)$.
 
 The construction in Lemma \ref{t2}  can be  applied to an arbitrary bijection
$T:\overline n \to \overline n$, not necessarily a chart.  The difference is that a vertex   $v\in \Gamma$ is now thickened to a union of
several disks meeting each other at the center; the number of these disks is equal to the number of $T$-orbits  contained in the set
$J=J(v)$. This yields a curve  with corners  on a {\it singular surface} $\Sigma^s=\Sigma^s(T)$ obtained from a closed (oriented, possibly
non-connected) surface
$\Sigma=
\Sigma (T)$ by contracting   several disjoint finite subsets.  The image of this curve a pre-clean dessin d'enfants on $\Sigma$ which
makes each component of
$\Sigma$  an algebraic curve over $\overline {\QQ}$
and  yields  a  point of the  corresponding moduli space. Consider in more detail  the case  where all singularities of  $\Sigma^s$ are {\it
simple} in the sense that 
$\Sigma^s$ is obtained from
$\Sigma$ by contracting  disjoint  pairs of   points.  It  happens iff  for any $T$-orbit $A\subset \overline n$ there is a $T$-orbit
$A'\subset \overline n$  such that  
$A\cup A'$ is invariant under negation on $\overline n$.  Under  this condition, $T$ gives rise to a modular graph $(\tau, g)$ in  the sense
of
\cite{ma}, p.\ 88. Its vertices are numerated by the components of $\Sigma$; its tails are numerated by negation invariant $T$-orbits;
its edges are numerated by unordered pairs $(A,A')$ of  $T$-orbits such that $A,A'$ are not negation  invariant but their union is. The
function
$g$ assigns to each vertex of $\tau$ the genus of the corresponding component of $\Sigma$.  In terminology of 
\cite{ma}, $\Sigma^s$   is the combinatorial type of the  \lq\lq prestable labeled curve" represented by $(\tau, g)$.
If the modular graph $(\tau, g)$ is stable, then our constructions give a   point  of the Mumford-Deligne compactification
$\overline {\mathcal M}_{G,k}$ where $k$ is the number of tails of $(\tau, g)$ and $G$ is the sum of all values of $g$ plus the first
Betti number of the 1-dimensional CW-complex underlying the graph $\tau$.

\section{Automorphisms}

\subsection{Automorphisms   of charts}\label{automo} An {\it automorphism}  of   a chart
$(n,t)$  is a permutation  $\varphi:\overline n \to
\overline n$ such that $\varphi t= t \varphi$ and $\varphi$ is a    power of the circular permutation
$ \sigma_n:\overline n \to
\overline n$.  The automorphisms of $(n, t)$ form a group with respect to composition, it is denoted $\Aut (t)$. This
group is a subgroup of the cyclic group of order $n$ generated by $\sigma_n$. Therefore $\Aut (t)$ is a cyclic group of
order dividing $n$.  Clearly,   $\Aut ( t)=\Aut(\sigma_n t (\sigma_n)^{-1})$. 

For example, consider the  charts  $t_1, t_2, t_3 :\overline 4 \to \overline 4$
given as products of two cycles
$$t_1=(1,3,-3,-1) (2,4,-4,-2),\;\;\;\;  t_2=(1,-1,3,-3) (2,-4,4,-2),\;\;\;\;  t_3=(1,-1,3,-3) (2,-2, 4, -4).$$
 It is easy to check that   $\Aut (t_1)= 1, \Aut (t_2)=\ZZ/2\ZZ$,
and
$\Aut (t_3)=\ZZ/4\ZZ$.

\begin{lemma} \label{ledivisibil}   For    any chart $(n,t )$ and   $m\geq
1$,  the order of    $  \Aut(t) $ divides $m\, k_m$ where $k_m$ is the number
of orbits of $t$ consisting of 
$2m$ elements.
   \end{lemma} 
\begin{proof}   Let $\overline M\subset \overline n$ be the union    of all orbits of
$t$ consisting of  
$2m$ elements.  Set $M=\overline M\cap \hat n$ where $\hat n=\{1,2,...,n\}$. By the definition of a chart,
$\overline M= M\cup (-M)$. Any automorphism  $\varphi:\overline n \to
\overline n$ of $t$  fixes $\overline M$ setwise.   
Since $\varphi$  is a  power of  $\sigma_n$,  it   fixes $\hat n$ setwise.  Therefore  $\varphi(M)=M$. 
The mapping $\varphi\vert_M:M\to M$ preserves the cyclic order on  $M$  induced by the standard cyclic order  $1\prec 2\prec ... \prec
n\prec 1$ on
$\hat n$. Therefore $(\varphi\vert_M)^{\card (M)}=\id$.  A  power of  $\sigma_n$
having a fixed point is the identity. Hence $\varphi^{\card (M)} =\id$. Taking as
$\varphi$ a generator of  the cyclic group $\Aut(t)$, we obtain that $\vert \Aut(t)\vert$  divides $\card (M) =(1/2) \card
(\overline M)=mk_m$. 
 \end{proof} 

 \subsection{Automorphisms   of curves}\label{auto}  An {\it automorphism} of  a  curve  $f :S^1\to
\Sigma $ is  a homeomorphism of $f$ onto itself, that is a pair 
 $\varphi=( \varphi^{(1)} :S^1 \to S^1, \varphi^{(2)} :\Sigma \to \Sigma)$  of orientation  preserving  
  homeomorphisms  such that $f 
\varphi^{(1)}=\varphi^{(2)} f$ and
$\varphi^{(2)}(\cor (f))=\cor (f)$.  These conditions imply that  $\varphi^{(2)}$   preserves  $f(S^1),
\cross (f)$, and $\cor (f)$ setwise.  Automorphisms
$\varphi=(
\varphi^{(1)},
\varphi^{(2)} )$
 and $\psi=( \psi^{(1)}, \psi^{(2)} )$ of $f$   can be multiplied by  
$ \varphi \psi =(\varphi^{(1)} \psi^{(1)}, \varphi^{(2)} \psi^{(2)})$. With this
multiplication, the automorphisms of $f$  form a group.  Two automorphisms   of $f$ are {\it isotopic} if they can be included in a  
continuous 1-parameter   family of automorphisms of $f$.   Quotienting the group of automorphisms of $f$ 
by isotopy we obtain a group $\Auut  (f)$ whose elements are   isotopy classes of automorphisms of
$f$.  The
following theorem computes 
$\Auut  (f)$   in terms of the chart of
$f$.

	\begin{theor}\label{t389}   For any filling   curve $f$ with chart $(n,t)$, we have  $\Auut  (f)=\Aut (t)$.
   \end{theor} 
\begin{proof}  We begin with a  simple and  well known
geometric lemma. The second claim of this lemma is    due to 
 J. W.  Alexander  and H. Tietze.

 \begin{lemma} \label{lemop1}   Let $D=D^N$ be a   closed
$N$-dimensional ball. Then

(i) any  homeomorphism  $g:\partial D\to \partial D$ extends to  a    homeomorphism  $\tilde g:  D\to D$ such
that under continuous deformation of $g$ its extension $\tilde g$   deforms continuously and for $g={\id}_{\partial D} $ we have $ 
\tilde {g} =\id_{D}$;

(ii) any  homeomorphism  $ 
D\to D$ fixing    $\partial D$ pointwise  is isotopic to the identity in the class of  homeomorphisms $D\to
D$ fixing   $\partial D$ pointwise.
   \end{lemma} 
\begin{proof}  We   identify $D$ with the unit Euclidean
ball   
 $ \{z\in {\RR^N}, \vert \vert z\vert\vert \leq 1\}$ where $\vert\vert.\vert \vert$ is the  Euclidean
norm on
$\RR^N$. Then $ \partial D= \{z\in {\RR^N}, \vert \vert z\vert\vert =  1\}$.  For a homeomorphism $g:\partial D\to \partial D$ and $
z\in D - \{0\}$, set $\tilde g (z) = \vert \vert z\vert\vert\, g(z/\vert \vert z\vert\vert)$. Set $\tilde g(0)=0$.  This 
defines a homeomorphism $\tilde g: D \to D$ extending $g$ and satisfying  (i).

For a homeomorphism  $h:
D \to D $ fixing    $\partial D$ pointwise, the formula
$$h_s (z) =
\left\{\begin{array}{ll}
z,  & \mbox{if\;\;   $s\leq  \vert \vert z\vert\vert \leq 1$}, \\
s\, h(z/s), & \mbox{if\;\;    $ \vert \vert z\vert\vert < s$} 
\end{array} \right.
$$
defines   an isotopy (= a continuous family of homeomorphisms) $\{h_s:D\to D\}_{s\in [0,1]}$ of $h_0=\id$ to $h_1=h$.
\end{proof}

Let $f$ be a	 filling   curve  with chart $(n,t)$. If $n=0$, then  $f$ is a trivial curve. Using Lemma
\ref{lemop1} (ii) and  the fact that   orientation preserving homeomorphisms $S^1\to S^1$ are isotopic to the identity, we  
obtain    $\Auut  (f)=1=\Auut  (t)$. Assume from now on that $n> 0$.  As we know,  any  automorphism $\varphi$   of    $f$
induces a   permutation
$\varphi_*:\Fl(f )
\to
\Fl(f )$.   Isotopic automorphisms of $f$ induce the same permutations. (Indeed, 
 a continuous deformation of a  permutation of $\Fl (f)$ must be constant since $\Fl (f)$ is finite.) Therefore the  formula $F(\varphi)
=\varphi_*$  defines a group homomorphism   $F:\Auut  (f)\to \Iso (\Fl(f))$ where for a set $S$, we denote by $\Iso (S)$
 the group of 
      bijections $S\to S$ with multiplication given by composition.

 \begin{lemma} \label{le1bnb}    The homomorphism $F: \Auut  (f) \to  \Iso (\Fl(f))$ is injective.  
   \end{lemma} 
\begin{proof}   Let  $\varphi=(\varphi^{(1)}, \varphi^{(2)})$ be an automorphism of $f$ whose isotopy class lies in the kernel of $F$.  We
shall show that
$\varphi$ is isotopic to the identity.  We prove  first  that $\varphi$ is isotopic to an automorphism $\eta=(\eta^{(1)}, \eta^{(2)})$ of
$f$ with
$\eta^{(1)}=\id_{S^1}$.  Since
$\varphi_*=F(\varphi)$  acts as the identity on all flags of $f$, the homeomorphism  $\varphi^{(1)}:S^1\to S^1$ must 
fix
$\Cross (f)$ pointwise.   The points of $\Cross (f)$ split $S^1$ into $n=n(f)$ consecutive arcs
$\gamma_1,...,
\gamma_n$.  Since   $\varphi^{(1)}$ is orientation preserving and
constant on   $\partial \gamma_k$,  it   maps  $\gamma_k$ onto itself for  all $k$.  By Lemma \ref{lemop1} (ii),  the
restriction of  $\varphi^{(1)}$ to $\gamma_k$     is isotopic to the identity in the class of
homeomorphisms  $\gamma_k \to \gamma_k$   fixing $\partial \gamma_k$ pointwise.  Gluing these isotopies over
$k=1,...,n$, we obtain a  continuous  family of homeomorphisms
$\{\varphi^{(1)}_s:S^1\to S^1\}_{s\in [0,1]}$ fixing  $\Cross (f)$ pointwise and such that 
$\varphi^{(1)}_0=\varphi^{(1)}, \varphi^{(1)}_1=\id$.  The graph    $f(S^1)$ is obtained from $S^1$ by
identifying  the points of $\Cross (f)$ having  the same image under $f$.    Since
the homeomorphism $\varphi^{(1)}_s:S^1\to S^1$ fixes
$\Cross (f)$ pointwise, it induces a homeomorphism $\psi_s: f(S^1)\to f(S^1)$ such that
$\psi_s f =f \varphi^{(1)}_s $,   $\psi_s$  fixes 
$\cro (f)$ pointwise,  and $\psi_1=\id$. Since $f$ is  a filling curve, the splitting  of $\Sigma$ along $f(S^1)$  gives  a finite
family  of closed 2-disks $\{D_q\}_q$.  Each    $ \psi_s  $ induces a  homeomorphism  $\partial D_q \to \partial D_q$.  By
Lemma
\ref{lemop1} (i), the latter extends 
  to a   homeomorphism  $
D_q \to  D_q$. Gluing these homeomorphisms along $f(S^1)$, we
obtain     a continuous  family of homeomorphisms $\{\tilde \psi_s: \Sigma\to
\Sigma\}_{s\in [0,1]}$ such that   $\tilde \psi_s \vert_{f(S^1)} =\psi_s $  for all $s\in [0,1]$  and $\tilde
\psi_1=\id_\Sigma$. Then
 $\tilde \psi_s  f =\psi_s f=f \varphi^{(1)}_s:S^1\to \Sigma$. Thus we obtained  a continuous  family
$(\varphi^{(1)}_s:S^1\to S^1,  \tilde \psi_s: \Sigma\to \Sigma)$  of automorphisms of $f$ such that
 $\varphi^{(1)}_0=\varphi^{(1)}$ and  $ (\varphi^{(1)}_1, \tilde \psi_1)=\id$.  The family 
$\{ \eta_s= (\varphi^{(1)} (\varphi^{(1)}_s)^{-1} ,  \varphi^{(2)} (\tilde \psi_s)^{-1} )\}_{s\in [0,1]}$ of automorphisms of $f$ is
an isotopy between  $\eta_0=(\id_{S^1},  \varphi^{(2)} (\tilde \psi_0)^{-1})$ and 
$\eta_1=\varphi$.  Thus $\varphi $ is isotopic to an automorphism $\eta=\eta_0$ of $f$  with 
$\eta^{(1)}=\id_{S^1}$. 
 
  We prove  now  that any such   $\eta $  
 is isotopic to the identity. The equality  
$ \eta^{(2)} f=f 
\eta^{(1)}=  f$ shows that $\eta^{(2)}:\Sigma\to \Sigma$ preserves  $f(S^1)$ pointwise. 
Splitting $\Sigma$ along $f(S^1)$   we 
obtain a finite family  of closed 2-disks $\{D_q\}_q$.     Since     $\eta^{(2)}   $  preserves
orientation, it has to map each  $D_q$  into itself fixing $\partial D_q$ pointwise. By Lemma \ref{lemop1} (ii), 
the resulting self-homeomorphism of 
  $D_q$  is isotopic to the identity in the class of self-homeomorphism of 
   $D_q$
 fixing $\partial D_q$  pointwise. Gluing such isotopies over all $q$, we obtain an isotopy of $\eta^{(2)}$ to the identity
  in the class of homeomorphisms $\Sigma\to \Sigma$ fixing $f(S^1)$ pointwise. All these homeomorphisms are
automorphisms of $f$ acting on $S^1$ as the identity.  Hence  $\eta $  
 is isotopic to the identity automorphism of $f$. 
\end{proof}

 \begin{lemma} \label{le99}   Under the    identification $\Fl (f)=\overline n$, the image of    $F:\Auut  (f) \to    \Iso (\Fl(f))= \Iso
(\overline n)$  is   $\Aut(t)$.
   \end{lemma} 
\begin{proof}    Pick   $\varphi\in \Auut  (f)$.  By  Section  \ref{action}, the bijection 
$F(\varphi)= \varphi_*: 
 \Fl (f) \to \Fl (f)  $ commutes with the flag rotation.  This bijection is defined   by    $\varphi_*((x, \pm))=
(\varphi^{(1)} (x),
\pm)$ where 
$x\in
\Cross (f)$.  The homeomorphism    
$\varphi^{(1)}:S^1
\to S^1$ sends  
$\Cross (f)$ onto
itself preserving the cyclic order.  Hence under the identification  $\Fl (f)=\overline n$, the bijection $\varphi_*$ becomes a power
of the circular permutation
$\sigma=\sigma_n:\overline n \to \overline n$. Thus, $\varphi_*\in
\Aut (t)$. 

It remains to prove that all the powers of  $\sigma $  commuting with $t$  lie in  the
image of $F$.     
Note that if  this  holds for  a  curve homeomorphic to
$f$, then this holds for $f$. Therefore it is enough to prove our claim  in the case where $f$ is the filling  curve
constructed from the chart $(n, t )$ in the proof of Lemma \ref{t2}.

We  use notation introduced there: the points 
$x_k=\exp (2k \pi i /n)\in
S^1$,  the graph $\Gamma $, and the
projection $p:S^1\to \Gamma$.  Let    $\tau=\sigma^q $ with $q\in \ZZ$  be a   power  of 
$\sigma $  commuting with
$t$.
Let
$\varphi^{(1)}:S^1\to S^1$  be  multiplication by $\exp (2q \pi i /n)$.   Clearly, $\varphi^{(1)}
(x_k) =x_{\tau(k)}$ for all $k$.  If    $p(x_k) =p(x_j)$ with  $k,j\in \hat n$,  then  
there is a power of
$t$ transforming $j$ into $k$.  Since $t \tau=\tau t$, the same power of $t$  transforms  ${\tau(j)}$ into
${\tau(k)}$.  Hence
 $p \varphi^{(1)} (x_k)=p(x_{\tau(k)})= p(x_{\tau(j)})=p \varphi^{(1)} (x_j)$. 
Therefore  $\varphi^{(1)} $ induces a  map 
$\psi:
\Gamma\to
\Gamma
$ such that
$\psi p= p\varphi^{(1)}$.   Applying the same argument to  the inverses of $\tau$ and $ \varphi^{(1)}$, we obtain 
an inverse  map. Hence $\psi$ is a homeomorphism. 

For a  vertex $v=p(x_k)\in \Gamma$,  the point  $\psi(v)= p\varphi^{(1)}(x_k) =p(x_{\tau(k)})$    is also a vertex of
$\Gamma$.   Clearly, 
$\psi$ maps the neighborhood
$V_v\subset
\Gamma$ of $v$ defined   in  Lemma \ref{t2}   onto  $V_{\psi(v)}$.   
The    homeomorphism $\psi\vert_{V_v}: V_v\to V_{\psi (v)}$ transforms  a  flag $(x_k,\pm)$ at $v$  into the
flag $(x_{\tau(k)}, \pm)$ at $\psi(v)$.  Since $t \tau=\tau t$,  this transformation of flags is $t$-equivariant. 
Therefore  $\psi\vert_{V_v}$ maps the arcs
forming
$V_v$ onto the arcs forming 
$V_{\psi (v)}$ preserving their cyclic order induced by $t$.  This implies that  
$\psi\vert_{V_v} $ extends to an orientation preserving     homeomorphism   $D(v) \to D
(\psi(v))$.  The resulting self-homeomorphism of $\Gamma\cup \cup_v D(v)$ obviously extends to the  ribbons used in   
Lemma \ref{t2} to construct  
$U$.  This gives   an orientation preserving  homeomorphism    $ U\to 
U$.  By
Lemma \ref{lemop1} (i), the latter   extends to a  homeomorphism  $\varphi^{(2)}:\Sigma \to \Sigma $.
  Since  $ \varphi^{(2)}$ is an extension of $\psi  $, we have  $i \psi =\varphi^{(2)} i $ where  $i$ is the inclusion
$\Gamma
\hookrightarrow \Sigma$. Recall that
$f=ip:S^1\to
\Sigma$.  Therefore $f \varphi^{(1)}=i p \varphi^{(1)}= i \psi p= \varphi^{(2)} i  p= \varphi^{(2)} f$. Hence  
$\varphi=(\varphi^{(1)}, \varphi^{(2)})$  is an automorphism of $f$. It follows from the definitions that
$F(\varphi)=\tau$. \end{proof} 

Lemmas \ref{le1bnb} and \ref{le99}  directly imply  Theorem \ref{t389}.
\end{proof}

\begin{corol}\label{cordivi}  For a filling    curve  $f$,  the group 
 $  \Auut  (f) $ is cyclic of finite order dividing   $m\,k_m(f) $ for all $m\geq 1$ where     $k_1(f), k_2 (f),...$  are the numbers defined
in Section \ref{regu}.
   \end{corol}

 \section{Self-transversal  curves and semicharts}

 \subsection{Self-transversal  curves}\label{self} A   curve
$f:S^1\to \Sigma$  is  {\it self-transversal} at   $a\in \cross (f)$  if for any  distinct   $x,y\in f^{-1} (a) \subset
S^1$   the branches of
$f$ at
$x$ and
$y$ are topologically transversal.  The latter condition means that there  is  a homeomorphism  of  a   neighborhood  of  
  $a$
in $ 
\Sigma$  onto   $   \RR^2$ whose composition with $f$ sends a neighborhood of $x$ (resp. of $y$) in $S^1$  to
$\RR\times 0$ (resp.  to 
$0\times
\RR$).    A curve   
$f $  is  {\it self-transversal} if it is self-transversal at all points of  $ \cro (f)$.

We now describe the charts of self-transversal  curves.  Let us say that a chart $(n,t)$ is {\it
straight} if
$t(-k)=-t(k)$ for all
$k\in
\hat  n$. This property   is preserved uner conjugation of $t$ by the circular permutation. 

	\begin{lemma}\label{lepsikk}  A   curve   is self-transversal if and only if its   chart is straight.
 \end{lemma}
\begin{proof}   For   a flag
$r\in \Fl(f)$ of a curve $f$, denote the opposite flag by $-r$. Under the identification $\Fl(f)=\overline {n}$ where $n=n(f)$  the
involution 
$r\mapsto -r$ on $\Fl(f)$ corresponds to the
negation on $\overline {n}$.  Therefore it suffices to prove that  $f$  is self-transversal if and only if 
$t(-r)=-t(r)$ for all  $r\in \Fl (f)$.  

   Moving  around a point   $a\in \cross (f)$    in the positive direction,  we  cyclically numerate the  flags
of $f$   at  this point  $r_1,r_2,...,r_{2m}$ where $ m \geq 2$ is the multiplicity of
$a$.   For every $q\in \ZZ/2m\ZZ$,   there is a unique $q' \in \ZZ/2m\ZZ$ such that 
$-r_q=r_{q'} $.    The curve 
$f$ is self-transversal  at $a$ if and only if    $q'={q+m}$ for all $q$.
  If $f$ is self-transversal at $a$, then 
 $$ t(-r_q)=t(r_{q+m})=r_{q+m+1}=-r_{q+1}=-t(r_q).$$  Conversely, suppose that $t(-r_q)=-t(r_q)$ for all  $q$. 
Then     $$r_{2q'-q}=t^{q'-q} (r_{q'})=t^{q'-q}    (-r_q) 
=-t^{q'-q} (r_q) =-r_{q'}=r_q.$$
Therefore $2(q'-q)=0 $.  If $q'=q$, then   $-r_q=r_{q'}=r_q$
which is impossible. Hence
$q'=  q+m$.   
   \end{proof} 

By  Theorem \ref{t5} and Lemma \ref{lepsikk},
the formula $f\mapsto (n(f),t(f))$  defines a bijective correspondence between    self-transversal   filling   curves 
considered up to homeomorphism  and  straight charts  
  considered up to conjugation by the   circular permutation.

 \subsection{Semicharts}\label{restc} Information contained
in a straight chart can be packaged in   a  more compact way. A {\it semichart} a triple $(n,v:\hat n\to \hat n,S\subset \hat n)$
where $n\geq 1$, $v$ is a bijection,
and $S$ meets  each orbit of $v$ in an odd number of elements. The   group $\ZZ/n\ZZ$ acts on the set of semicharts  
via    $(n,v,S)\mapsto (n, \sigma v \sigma^{-1}, \sigma (A))$ where
$\sigma=\sigma_n\vert_{\hat n}:\hat n\to
\hat n$. 
  By convention, there is   a unique   {\it empty semichart} $(n,v,S)$ with $n=0$.

	\begin{lemma}\label{lepsi} For each $n\geq 0$, there is a $(\ZZ/n\ZZ)$-equivariant bijective correspondence between
straight charts
$(n,t)$ and semicharts $(n,v, S)$.  Under this correspondence $\overline n/t=\hat n/v$.
 \end{lemma}
\begin{proof}  It suffices to consider the case $n\geq 1$.  With a straight chart $(n,t)$ we associate a semichart
  as follows. For $k\in \hat n$, set
 $v (k) = \vert t(k) \vert \in \hat n$.  
Set   $S=\{k\in \hat n\,\vert \,  t(k)<0\}$.  We check that $(n,v,S)$ is a semichart.   If $v(k)=v(l)$ for $k,l\in
\hat n$, then   either $t(k)=t(l)$ or $t(k)=-t(l)=t(-l)$. Since $t$ is a bijection,  
$k=l$ or 
$k=-l$. The latter is impossible since $k,l>0$. Thus   $v$ is injective and therefore bijective.  Let 
$A\subset \hat n$ be an orbit  of
$v$.   Set
$m=\card (A)$ and pick  $a\in A$.  Then $v^m(a)=a$ and $A=\{a, v(a),...,
v^{m-1} (a)\}$. An   induction on $q=1,2,...$ shows
that    $t^q(a)=(-1)^w v^q(a)$ where  $w$ is the number of terms of the sequence $a, v(a),...,
v^{q-1} (a)$ belonging to $S$. 
Therefore  $t^m(a)=(-1)^w v^m(a)=(-1)^w a$ where  $w=\card (A\cap S)$.  The
equality $t^m(a)=a$ would imply that the orbit of $t $  containing
$a$ is contained in the set  $\{a,t(a),..., t^{m-1} (a)\}$. All  its  elements are distinct from $-a$ since the  absolute values 
of $t(a),..., t^{m-1} (a)$ are equal to 
$v(a),..., v^{m-1} (a)$, all distinct from $a$. This   contradicts the definition of a chart. Hence
$t^m(a)=-a$ so that   $w=\card (A\cap S)$   is  odd.  It is clear that the orbits of $v$ are
obtained by intersecting the orbits of $t$ with $\hat n$. This gives a bijection 
$\overline n/t=\hat n/v$. 

Conversely, having a semichart $(n,v,S)$ we define a map $t:\overline n\to \overline n$ by  
$t(k)=v(k)$ for   $ k\in \hat n - S$,   $t(k) = -v(k) $ for  $k\in S$, and   $t(k)=-t(-k)$ for $k\in -\hat n$. Then $t(\overline n)=\{\pm
v(k)\}_{k\in
\hat n}=\overline n$, since $v( \hat
n)=\hat n$. Thus $t$ is   bijective. The arguments above show that for any $a\in \hat n$, we have $t^m(a)=-a$ where $m$ is the
number of elements in the $v$-orbit of $a$. Hence $(n,t)$ is a chart.  By its very definition, it is straight. It is   clear that the arrows
from charts to semicharts and backward  defined above are mutually inverse.  The equivariance  of these arrows with respect to the
conjugation by the circular permutation  is straightforward.     
\end{proof} 

Example:  the  chart   $ \overline 4 \to \overline 4$
defined  as the product  of two cycles
 $ (1,3,-1,-3) (2,4,-2,-4)$ is straight. The associated semichart $(4, v , S)$ is  
$v=(13) (24)$, $S=\{3,4\}$.

Using Lemma \ref{lepsi}, we can associate a semichart (considered up to conjugation by   the circular permutation) with any
self-transversal  curve. Combining with the results of Section
\ref{self},  we obtain the following.

\begin{theor}\label{t12658745555} There is a bijective correspondence between  self-transversal   filling   curves 
considered up to homeomorphism  and  semicharts  
  considered up to conjugation by   the   circular permutation.
   \end{theor}

 \subsection{Automorphisms}\label{autosemic}    An {\it automorphism}  of   a    semichart
$(n,v,S)$ with $n\geq 1$  is a permutation  $\mu:  \hat n \to
\hat n$ such that $\mu v= v \mu$, $\mu(S)=S$, and $\mu$ is a  power of  
$ \sigma_n \vert_{\hat n}:\hat n \to
\hat n$.  The automorphisms of $(n, v, S)$ form a group with respect to composition denoted $\Aut (v,S)$.  By convention, for  the 
empty semichart, this group is trivial.   If  $(n,t)$  is the straight chart determined by  $(n,v,S)$, then restricting   automorphisms   of
$(n,t)$ to
$\hat n$, we obtain  that 
$\Aut(t)=\Aut (v,S)$.
This and Theorem \ref{t389} give the following.

	\begin{corol}\label{t3892nnn}   For a    self-transversal  filling    curve $f$ with semichart $(n,v,S)$, we have  $\Auut  (f)=\Aut (v,S)$.
   \end{corol}

Analyzing the automorphisms of  a semichart $(n,v,S)$,    it is easy to see that an   automorphism    preserving an orbit of $v$
setwise   has an odd order.  We prove the corresponding fact for curves.
 
   \begin{theor} \label{le356g4}    Let $f$ be a  self-transversal   filling 
 curve and  $a\in \cro (f)$.  Let $\Aut_a(f)$  be the subgroup of $\Aut(f)$ consisting of the (isotopy
classes of) automorphisms of $f$ preserving $a$.  Then $\Aut_a(f)$ is a finite cyclic group of odd order dividing 
  the multiplicity $m=m_a$
  of
$a$.  
   \end{theor} 
\begin{proof}   Moving around    $a $ in the positive direction,  we   cyclically numerate $r_1,r_2,...,r_{2m}$ the  flags of $f$  
at  $a$.   Let $j$ be the minimal  element of the set $\{1,2,..., 2m\} $ such that there is $\varphi\in
\Aut_a(f)$ with $\varphi_*(r_1)= r_{1+j}$.   Then $\varphi_* (r_q)=r_{q+j}$ for all $q\in \ZZ/2m\ZZ$ and  any  $\psi \in \Aut_a(f)$ is
a power of $\varphi$. Indeed,   
$\psi_*(r_1)=r_q$ for some  $q$.   If  $q\neq 1 (\modu j)$ then a product of $\psi_*$ with a power of
$\varphi_*$ transforms $r_1$ into
$r_{1+k}$ with $1\leq k< j$ which contradicts the choice of $j$.  If $q={1+kj}$ with $ k\in \ZZ$, then  $\varphi_*^{-k}
\psi_*(r_1)=r_1$ and therefore $
\psi=\varphi^{k}$.  

By the choice of $j$, the residues
$ 1+2j, 1+3j,... \,(\modu 2m) $ do not take   values  $2,3,...,j-1$. This is possible only if  $dj=0 (\modu 2m)$ for some
$d\geq 2$.  Take the smallest such $d$. Then 
$\varphi_*^d (r_1)=r_{1+d j}=r_1$.  Hence  $\varphi^d=1$. It remains to prove that $d$ is odd. If $d=2e$ with
$e\in \ZZ$, then  $ej=m (\modu 2m)$ and  $\varphi_*^e (r_1)=r_{1+ej}=r_{1+m} $. Since $f$ is self-transversal,  the  flag 
$r_{1+m}$ is opposite to  $r_1$.   This  contradicts the equality 
$\varphi_*^{e} (r_1) = r_{1+m} $ since   automorphisms of curves cannot transform incoming  flags into outgoing ones.
 \end{proof}

 \begin{remar}\label{rista}  It is easy to see  that for any   curve $f$ and any  $a\in \cro (f)$ of multiplicity
$m \geq 1$,  we have 
$ \vert\Aut(f)/\Aut_a(f)
\vert\leq k_m(f)$.   
By Theorem \ref{le356g4},  if $f$ is self-transversal and filling, then      $\vert \Aut_a(f) \vert$ is an odd divisor of    $m$. 
In particular,  if  $m$ is a power of $2$, then  $\Aut_a(f) =1$. 
Under further assumptions on $f$, these observations may  give $\Aut(f)=1$.  For example, if  $k_{2^q}(f)=1$ for some $q\in \ZZ$, then
$\Aut(f)=1$.  Another example: if
$k_5(f)=4$ and
$k_6(f)=1$, then    Corollary \ref{cordivi} and Theorem \ref{le356g4}  imply that $\Aut(f)=1$. 
\end{remar}

 \subsection{Coxeter groups}\label{coxx} The permutations $t:\overline n \to \overline n$ such that $t(-k)=-t(k)$ for all
$k\in
\overline n$ form a group  $W_n$. This is   the Coxeter group  of type $B$, see for instance \cite{bo}, \cite{bb}.  By the results
above,  the charts of   self-transversal curves  are elements of $W_n$ for an appropriate $n$. (To avoid the circular
indeterminacy in the definition of the charts, we  assume  the curves to be pointed, cf. Section  \ref{poin}.)
Conversely, every  element of $W_n$ whose orbits in $\overline n$ are negation invariant  gives rise  to a pointed self-transversal filling
curve on a surface of a certain genus.  Other  elements of $W_n$ give rise to  pointed self-transversal filling curves on  singular surfaces
with simple singularities, cf. Section
\ref{encoreeee}.    This  somewhat surprising connexion  between curves and Coxeter groups of type $B$  can be exploited to study  
curves. For example, for any pair of  pointed self-transversal  curves (possibly lying on different surfaces but yielding  the same   number
$n$),  we can consider the Kazhdan-Lusztig polynomial  of their charts. It seems however that it is quite difficult to compute
this polynomial  directly from the curves. This is due to a poor   connexion  between  multiplication  in $W_n$ 
and the topology of   curves.

 \section{Counting curves}\label{countdooo}

\subsection{Preliminaries} We begin with a few  simple remarks on   group actions.  The set of orbits of a (left) 
action of  a group $G$ on a set $\mathcal S$ is denoted $\mathcal S/G$. For $a\in \mathcal S$, let 
 $\Stab (a)=\{g\in G\,\vert\, ga=a\}$ be the stabilizer of  $a $. For an  orbit
 of this action $A \subset \mathcal S$,   set  $\Stab_A = \Stab (a)  $ for some  $a\in A$.  The  isomorphism class of the
group 
$\Stab_A
 $ does not depend on the  choice  of $a$:   if  $a'\in A$, then  there is $g\in G$ such that $a'=ga $ and  $\Stab (a')
=g\Stab(a) g^{-1}$. The number of elements of a finite group $G$ will be  denoted
$\vert G\vert$.

\begin{lemma} \label{le5263}    Let $G$ be a finite group acting  on a finite set $\mathcal S$. 
Then
\begin{equation}\label{osno}\card (\mathcal S)=\vert G\vert \sum_{A\in \mathcal S/G} \frac{1} {\vert \Stab_A\vert} .\end{equation} \end{lemma} 
\begin{proof}   For the orbit   $A\subset \mathcal S$ of   $a\in \mathcal S$, the formula $g\mapsto ga$ defines a  surjection
$G\to A$   and a bijection $G/ \Stab (a)\approx A $. Hence  $\card (A)= \vert G\vert / \vert  \Stab
(a)\vert=\vert G\vert/{\vert \Stab_A\vert}$.  Summing up over all the orbits, we obtain  $$\card (\mathcal S)=\sum_{A\in \mathcal S/G} \card (A)   
=\sum_{A\in \mathcal S/G}
 \vert G\vert/{\vert \Stab_A\vert}=\vert G\vert \sum_{A\in \mathcal S/G}
 1/{\vert \Stab_A\vert} .$$ \end{proof}

\subsection{Counting  filling curves}  A sequence of integers $K=(k_1,k_2,...)$ is {\it finite} if $k_m=0$ for all
sufficiently big $m$. For a finite sequence $K$, set  
$n(K)=\sum_{m\geq 1} m k_m$.  By convention,   $(-1)!=0!=1$. 

\begin{theor}\label{thefo}  Let   $K=(k_1,k_2,...)$ be a finite sequence  of non-negative integers. Let $\mathcal C(K)$
be the set of  homeomorphism classes  of  filling  curves $f$
such that    $k_m(f)=k_m$  for  all $m\geq 1$. Then
\begin{equation}\label{eefroo}\sum_{f\in \mathcal C(K)} \frac {1} {\vert \Auut  (f)\vert}  
=  (n(K)-1)!  \prod_{m\geq 1} \frac{1} {k_m!}  \left ( \frac
{(2m-1)!}{    m!}\right )^{k_m}   .\end{equation}
   \end{theor} 
\begin{proof}  If $k_1=k_2=...=0$, then both sides of this formula are equal to 1.  Assume   that
at  least one   $k_m$ is non-zero so that $n=n(K)\geq 1$.  Let
$\mathcal S
$ be the set of  charts
$(n, t)$ such that 
  $t$   has    $k_m$
orbits   of cardinality
$2m$ for all
$m\geq 1$.   The set $\mathcal S$ is invariant under the   action of  
$G=\ZZ/n\ZZ$ on charts defined in Section \ref{charts}.   
By Theorem \ref{t5},  assigning to   $f\in \mathcal C(K)$ its chart,   we obtain  
$\mathcal C(K)=\mathcal S/G$. By Theorem \ref{t389},  the stabilizer of the chart of $f$ is isomorphic to  $ \Auut
(f)$. Hence Formula \ref{osno} gives
\begin{equation}\label{eefr1}\card (\mathcal S)=n \sum_{f\in \mathcal C(K)} \frac{1} {\vert \Auut (f) \vert} .\end{equation}
We now compute $\card (\mathcal S)$. To specify $t\in \mathcal S$ we   need to  specify a partition  of $\overline n$ into the orbits of
$t$ and   the action of $t$ on these orbits. Since the  orbits of $t$ are  invariant under the negation, the  partitions of  $\overline n$ in
question bijectively correspond  to   splittings of the set
$\hat n=\{1,2,...,n\}$ into $k_1+k_2+...$ disjoint subsets such that $k_1$ of them have 1 element, $k_2$ of them have 2
elements, etc.    The number of such splittings of  $\hat n$ is
\begin{equation}\label{eefrici1} \frac {n!}{\prod_{m\geq 1} (m!)^{k_m} k_m!}.\end{equation}
The number of 
transitive actions of 
$t$ on   a set  of $M\geq 2$ elements    is   $(M-1)!$.
Therefore $$\card (\mathcal S)= \frac {n!}{\prod_{m\geq 1} (m!)^{k_m} k_m!} \prod_{m\geq 1} ((2m-1)!)^{k_m}
= {n!}  \prod_{m\geq 1} \frac{1}{k_m!  } \, ({(2m-1)!}/{    m!})^{k_m}  .$$
This equality  and  Formula \ref{eefr1} imply   Formula \ref{eefroo}.
  \end{proof} 

 \begin{corol} \label{codot1}    If under the conditions of Theorem \ref{thefo}, 
$gcd \{mk_m\}_{m\geq 1}=1$, then 
$$\card {  \mathcal C (K) }  
= {(n(K)-1)!}  \prod_{m\geq 1} \frac{1}{ k_m!  }  \, ({(2m-1)!}/{    m!})^{k_m}   .$$
   \end{corol} 

This follows directly from  Corollary \ref{cordivi} and Theorem \ref{thefo}. As an illustration   consider a few  special cases.
 In the case where  $k_m=0$ for all $m\geq 3$, we have
$n=k_1+2k_2$ and   Theorem \ref{thefo} gives 
\begin{equation}\label{sert} \sum_{f\in \mathcal C(K)} \frac {1} {\vert \Auut  (f)\vert}  
=  \frac{(k_1+2k_2 -1)!}{k_1!\, k_2! } \,3^{k_2}.\end{equation}
For $k_1=0$,  the right-hand side here simplifies to $ ({(2k_2 -1)!}/{ k_2! })  \,3^{k_2}$. For  $k_2=0$, the right-hand side of
Formula
\ref{sert} simplifies to $1/k_1$.  The resulting formula
  can be verified directly since the set  $  \mathcal C (k_1,0,0,...)$ consists of   one element $f$
that is  an embedding
$S^1\hookrightarrow S^2$ with
$k_1$   corners and  $  \Auut  (f) =\ZZ/k_1\ZZ$.

Formula \ref{eefroo} can be rewritten as an equality   in the ring $\QQ[[t_1,t_2,...]]$ of formal power
series in commuting variables $t_1,t_2,...$ with rational coefficients. Let  $\mathcal C$ be  the set of 
homeomorphism classes of filling  curves. 
Then
\begin{equation}\label{formal1}\sum_{f\in \mathcal C} \, \frac { 1 } { (n(f)-1)!\,\vert \Auut  (f)\vert } \, t_1^{k_1(f)}
t_2^{k_2(f)} t_3^{k_3(f)}\cdots   =       \prod_{m\geq 1} 
\exp 
\left ( \frac {(2m-1)!}{    m!}\, t_m \right )    .\end{equation}
Here the product on the right hand side   is   a limit of finite products
$\prod_{m\geq 1}^N $ when $N\to \infty$. A  typical monomial $ t_1^{k_1 } t_2^{k_2  }
 \cdots t_q^{k_q}$ appears  in this limit   with the same coefficient 
as in $\prod_{m\geq 1}^q$.

\begin{theor}\label{thefo2}  Let   $K=(k_1,k_2,...)$ be a finite sequence  of non-negative integers. 
Let $\mathcal C_{str}(K)$
be the set of  homeomorphism classes  of  self-transversal  filling  curves $f$
such that    $k_m(f)=k_m$  for  all $m\geq 1$. Then
\begin{equation}\label{eefr}\sum_{f\in \mathcal C_{str}(K)} \frac {1} {\vert \Auut  (f)\vert}  
=   (n(K)-1)! \prod_{m\geq 1} \frac{1} { k_m!}\left  (\frac{ 2^{m-1}}{m}\right )^{k_m}   .\end{equation}
   \end{theor} 
\begin{proof} The proof   goes along the same lines as the proof of Theorem \ref{thefo} except that here we
  count     semicharts.   To specify a semichart   $(n,v,S)$ we   specify a partition of $\hat n $ into the orbits
of
$v$,   the action of $v$ on these orbits, and the intersections of $S$ with the orbits.
The number of partitions of $\hat n$   into the orbits of
$v$ is given by Formula
\ref{eefrici1}.    An   induction on $m$ shows that a set of $m$ elements  contains  $2^{m-1}$  subsets   having an  odd
number of elements. Therefore the number of semicharts $(n,v,S)$ is equal to 
$$ \frac {n!}{\prod_{m\geq 1} (m!)^{k_m} k_m!} \prod_{m\geq 1} (2^{m-1} (m-1)! )^{k_m}
=   n! \prod_{m\geq 1} \frac{1} { k_m!}\, (2^{m-1}/m)^{k_m}  .$$
 The rest of the argument is  as in the  proof of Theorem \ref{thefo}. \end{proof} 

 \begin{corol} \label{codot12}    If   $gcd\,
 (\{mk_m\}_{m\geq 1})=1$ of $k_{2^q}=1$ for some $q=0,1,2,...$, then 
$$\card {  \mathcal C_{str} (K) } 
=   (n(K)-1)! \prod_{m\geq 1} \frac{1} { k_m!}\, (2^{m-1}/m)^{k_m}
.$$
   \end{corol}

  Theorem \ref{thefo2} implies  Formulas \ref{formal2} and   \ref{eefr25}    of the introduction.

 \section{Further classes of curves}

\subsection{Pointed curves}\label{poin}  A {\it pointed curve} is a  curve $f $  
endowed with a distinguished point   $x\in S^1- \Cross (f)$  called the {\it base  point}. It should be stressed that   $x$ is
not viewed as a corner  of $f$. Using
$x$  as the starting point in the constructions of Section \ref{chartscurv}, we obtain a canonical identification $\Fl(f)=\overline n$ where
$n=n(f)=\card (\Cross (f))$. Therefore with a pointed curve   we can associate a  chart   without any indeterminacy.  A  {\it
homeomorphism} of pointed curves  is a homeomorphism of curves  mapping the base point   to the base
point. It is clear that homeomorphic pointed curves have  the same charts.  This yields    a bijective correspondence between 
pointed   filling   curves  considered up to homeomorphism  and   charts.
Similarly, there is   a bijective correspondence between 
 self-transversal pointed   filling   curves  considered up to homeomorphism  and  semicharts.
 
Note    that the group of  isotopy classes of  automorphisms    of a pointed filling curve is trivial. 

\subsection{Generic curves}\label{genc}   A  curve    is  {\it generic} if it is self-transversal and  all its crossings have multiplicity  
$2$.  We allow generic curves to have corners.       It is obvious that a  curve    is
generic if and only if  its chart $(n,t)$ is straight and  either $n=0$ or $n\geq 1$ and every orbit  of $t$  consists of   2 or 4
elements.   A self-transversal   curve    is generic if and only if  either $n=0$ or $n\geq 1$ and the map $v:\hat n\to \hat n$  in its
semichart 
$(n,v,S)$  is an involution.   We call such semicharts    {\it involutive}.  For an involutive semichart $(n,v,S)$, the   condition
that $S $ meets each orbit of   $v $  in an odd number of elements means  simply  that $S$ meets every orbit of $v$  in   one element.
Theorem \ref{t12658745555} yields  a bijective correspondence between  generic  filling   curves  considered up to homeomorphism  and 
involutive semicharts  
  considered up to conjugation by   the circular permutation.

\subsection{Alternating   curves}\label{altselftra}   A  curve  $f $  is  {\it alternating} if   the  
flag rotation $t:\Fl(f)\to \Fl (f)$ transforms  incoming flags into  outgoing  ones. Since the number of incoming and outgoing flags of $f$ 
is the same,
 $t$ then transforms outgoing flags into incoming ones.        A generic curve is alternating if and only if it has no crossings.  It is obvious
that a   curve  with chart
$(n,t)$   is alternating  if and only if   either  $n=0$ or $n\geq 1$ and $t(\hat n)=-\hat n$.

    A self-transversal curve $f$ with semichart $(n,v,S)$ is
alternating    if and only if   $n=0$ or 
$S=\hat n$.   The   condition that $S=\hat n$ meets each orbit of   $v $  in an odd number of
elements means simply that each orbit  of $v$ has an odd number of elements. In other words, 
$v$ must be  a permutation of odd order.  
Theorem \ref{t12658745555} yields  a bijective correspondence between  alternating self-transversal   filling   curves  considered up to
homeomorphism  and  permutations
$\{\hat n\to \hat n\}_{n\geq 0}$ of odd order
  considered up to conjugation by   the circular permutation. (By convention, for $n=0$, there is one permutation $\hat n\to \hat n$ of
odd order.)

\subsection{Beaming   curves}\label{becu} The notion of a beaming curve is in a sense opposite to the one of an alternating curve. 
A curve $f:S^1\to \Sigma$ is {\it beaming at a crossing} $a\in \cross (f) $ if one can draw a line on $\Sigma$ 
through
$a$   such that   the incoming flags of $f$ at   $a$ lie on one side of this line and   the outgoing flags of $f$ at  $a$ lie on
its other side.  We can rephrase this condition in terms of the flag rotation $t:\Fl(f) \to \Fl(f)$ by
saying that there is only one  incoming  flag $r$ at  $a$  such that $t(r)$ is outgoing (equivalently, 
there is only one  incoming  flag $r$ at  $a$  such that $t^{-1}(r)$ is outgoing).
The curve
$f$ is  {\it beaming}  if it is beaming at all its crossings.     An alternating curve is beaming if and
only if  it has no crossings. All generic curves are beaming.

 It is clear  that a   curve  with chart
$(n,t)$   is  beaming  if and only if   either $n=0$ or  $n\geq 1$ and    each orbit of   $t:\overline n\to \overline n$  contains only one 
element
$r>0$ such that  
$t(r)<0$.  
A self-transversal curve   with semichart $(n,v,S)$ is
beaming    if and only if   either $n=0$ or 
$S $  meets every orbit of   $v $  in one element.  

\subsection{Coherent   curves}\label{pecu} Consider a  self-transversal       curve $f $ and a crossing $a\in \cross (f)$.  A {\it branch  of
$f$ at $a$} is a union of two opposite flags   at $a$.  The  curve $f$ has $ m_a$
pairwise transversal branches    at
$a$.  Moving on the ambient surface in the positive direction around $a$ we obtain a cyclic order on the
set, $B_a$,  of branches of $f$ at $a$.  On the other hand, starting at a generic point  on the curve and traversing the whole  curve we   go
once along each branch of $f$ at $a$.    This also yields  a cyclic order on    $B_a$.    The curve 
$f$ is {\it coherent  at}  $a$ if these two cyclic orders coincide.  
A       curve  is {\it coherent} if it is  self-transversal   and   coherent at all    crossings.   Since a set of
two elements  has  only one cyclic order,   all generic curves   are coherent.  The  curve   in Example \ref{exam}.3 is coherent and
non-generic.     Note that if a coherent curve is non-generic, then inverting orientation on the curve (or on the ambient surface) we obtain
a non-coherent curve.

A semichart $(n,v,S)$ with $n\geq 1$  is {\it coherent} if  for any $k\in \hat n$ either $k<v(k)$ or  $k$ is the maximal element in its
$v$-orbit and then
$v(k)$ is the minimal element in the
$v$-orbit of $k$. By convention, the empty semichart is coherent. The coherency  of a
semichart  is preserved under conjugation  by the circular permutation.    All involutive semicharts are coherent.  

\begin{theor}\label{t5pontediksyyy}  A self-transversal     curve
  is coherent if and only if its semichart
   is coherent.
 This gives  a bijective correspondence between 
coherent     filling   curves  considered up to
homeomorphism    and  coherent semicharts  considered up to conjugation by   the circular permutation.
   \end{theor} 
\begin{proof} Let $f$ be a   self-transversal       curve with  chart $(n,t)$ and  semichart $(n,v,S)$.   Pick   $a\in \cross (f)$  and identify 
$\Fl (f)=\overline n$ as above. Let
$k_1<k_2<...<k_m$ be     elements of  $ 
\hat n$ corresponding to the outgoing flags at $a$.     The curve  $f$ is coherent at $a$ if and only if  $t(k_i)=\pm
k_{i+1}$ for
$i=1,...,m-1$ and $t(k_m)=\pm k_{1}$.  This holds iff   $v(k_i)= k_{i+1}$ for
$i=1,...,m-1$ and $v(k_m)= k_{1}$. Therefore $f$ is coherent at all its crossings iff  $(n,v,S)$ is coherent. The second claim of
the theorem follows    from the first claim and Theorem \ref{t12658745555}.
\end{proof}

\subsection{Perfect   curves}\label{pecui} A    pointed       curve $f$ is {\it perfect} if it is self-transversal, beaming, coherent, and 
satisfies the following property: $(\ast)$  starting at the base point   and moving along the curve we       enter each crossing
      $a\in \cross (f)$ for the first time along the unique incoming flag $s_a$   at $a$ such that the flag $t^{-1} (s_a)$ is outgoing. 
Under these assumptions,  
 the outgoing  flags of    $f$  at  $a$  can be numerated $r_1=-s_a ,r_2,..., r_m$ so that moving around $a$ we
encounter consecutively $r_1,r_2,..., r_m, -r_1, -r_2,..., -r_m$.  (Here    $t^{-1} (s_a)=r_m$.)    Note that if  
a curve   is   self-transversal, beaming, and coherent, then it is always possible to choose a base point  on it to satisfy $(\ast)$ at any
given crossing.  A perfect curve satisfies  $(\ast)$  for all crossings  and one and the same base point.  Trivial curves and trivial curves
with corners are perfect.  

A semichart $(n,v,S)$ is {\it perfect} if it  is coherent and  either $n=0$ or   
  $S $  meets every orbit of   $v:\hat n\to \hat n $ in one element which is the maximal element of this orbit with respect to the
standard order on
$\hat n\subset \RR$. 

\begin{theor}\label{t5pontediksyyyer}  A  pointed       curve
  is perfect iff its semichart
   is perfect.
 This gives  a bijective correspondence between 
perfect  pointed   filling   curves  considered up to homeomorphism  and  perfect semicharts.
   \end{theor} 
\begin{proof} The first claim     follows from the results above and the following  observation.  Consider  the chart $(n,t)$
of a pointed curve 
$f$ and   the
$t$-orbit $A\subset \overline n$ corresponding to
$a\in \cross (f)$.  Condition
$(\ast)$ can be interpreted by saying  that the only element $r\in A\cap \hat n$ such that $t(r)<0$ is the  maximal element of  $A\cap
\hat n$.    The second claim of the theorem follows  directly from the first claim.
\end{proof}

\subsection{Remark} The technique of charts yield  formulas counting   alternating/beaming/coherent/perfect
curves with weights as in Section \ref{countdooo}.  We state these formulas for  pointed curves.   Let $\mathcal C_p$
be the set of homeomorphism classes of pointed filling curves and $\mathcal C_{p, str} $ be  the set of homeomorphism classes of 
self-transversal pointed filling curves.   Then  $$\sum_{f\in \mathcal C_p} \, \frac { 1 } {  n(f) ! } \,  t_1^{k_1(f)} t_2^{k_2(f)} \cdots  
=     
\prod_{m\geq 1} 
\exp 
\left ( \frac {(2m-1)!}{    m!}\, t_m \right ),   $$
$$ \sum_{f\in \mathcal C_{p, str}} \, \frac { 1 } {  n(f) ! } \,  
t_1^{k_1(f)} t_2^{k_2(f)} \cdots   =      \prod_{m\geq 1} 
\exp 
\left ( \frac { 2^{m-1}}{    m}\, t_m \right )   .$$

 \section{Words}

\subsection{Words  and their automorphisms}\label{alphab}  An {\it alphabet}  is a finite set and   
{\it 
letters} are its elements.     A {\it  (signed)  word of length} $n\geq 1$  in an alphabet  $E$ is a pair (a mapping 
$w: \hat n\to E $, a set $S\subset \hat n$)  where $\hat
n=\{1,2,...,n\}$.    Such a 
word
 is encoded by the  sequence $w^S_1
w^S_2...w^S_n$ where
$   w^S_k= {w(k)} $ for $k\in S$ and  $ w^S_k= (w(k))^+$ for $k\in \hat n- S$. For example, the sequence
$AB^+A^+  $ in the alphabet $E=\{A,B\}$ encodes  the map  $ \hat 3 \to E$ sending $1, 2, 3 $
to
$A,B,A $ and  the set $S=\{1\} $.  More standard   unsigned words appear   when  $S=\hat n$.

Recall the map  $\sigma=\sigma_n:\hat n \to \hat n$
sending
$   k$ to $ k+1 $ for $k=1,2,..., n-1$ and sending
$   n$ to ~$ 1$. The {\it circular permutation} of  a word $(w, S)$ is  the
  word 
$(w \sigma^{-1} , \sigma  (S))$. This transforms   $w^S_1
w^S_2...w^S_n$ into  $w^S_n w^S_1
w^S_2...w^S_{n-1}$.  A 
word
$(w , S)$ is {\it full} if   $w$ is surjective, i.e.,  if  all letters appears in  
$w^S_1
w^S_2...w^S_n$  at least once (possibly with +).

The basic equivalence relation in the class of words is {\it congruence}. A word $(w:\hat n\to E, S)$
   is  {\it congruent} to a word $(w':\hat n'\to E', S')$     if  $n=n', S=S'$, and
there is a bijection $\psi:E\to E'$ such that $\psi w=w' $.  For  example, the  word $ABA^+$ in the alphabet $\{A,B\}$
is congruent to the   word $CDC^+$ in the alphabet $\{C,D\}$.
It is easy to classify   words up to congruence. Observe   that a  mapping $w:\hat n \to E$   gives rise to an equivalence relation
on   
$\hat n$ called {\it $w$-equivalence}:  
  $k,l\in \hat n$ are   $w$-equivalent if 
$w (k)=w(l)$.   If $w(\hat n)=E$ then the set of $w$-equivalence classes can be identified
with  
$E$ by assigning to a $w$-equivalence class its image in $ E$. It follows from the definitions that 
two full words  $(w , S)$, $(w' , S')$ of the same length $n$ are congruent if and only if $S=S'$ and the relations of  
$w$-equivalence and $w'$-equivalence on $\hat n$ coincide.

 An {\it automorphism}  of  a word $W=(w:\hat n\to E, S)$  
  is a  pair (a bijection   $\psi:E\to E$,  a residue $m\in \ZZ/n \ZZ$)  such that $\psi w=w \sigma^m:\hat n \to E $     and 
$  \sigma^m(S)=S$.  These  conditions   can
be rephrased by saying that the sequence   
$(\psi w )^{ S}_1
(\psi w )^{ S}_2...(\psi w )^{ S}_n$  can be obtained from  $w^S_1
w^S_2...w^S_n$   by the $(-m)$-th power of the circular permutation.
   The automorphisms
of $W$ form a group  $\Aut (W  )$  with unit $(\id, 0)$ and multiplication $(\psi, m) (\psi', m')= (\psi\psi', m+m')$.  
 This group    is preserved
under circular permutations  of $W$. The formula $(\psi, m) \mapsto m$ defines a \lq\lq forgetting"  group
homomorphism $p:\Aut(W) \to \ZZ/n\ZZ$.   

We   formally introduce a unique {\it empty word}   of length $0$ in an empty alphabet.  By convention, this word is full  
and its group of automorphisms  is trivial.

 \begin{lemma} \label{le1}   Let  $W=(w,S)$  be a  full   word of length $n\geq 1$ in an  alphabet $E$. Then

(i)  the  homomorphism 
$p:\Aut(W) \to \ZZ/n\ZZ$ is injective so that  $\Aut(W)$  is a finite cyclic group;

(ii)   the order of    $  \Aut(W) $ divides $gcd (\{m\, k_m\}_{m\geq 1})$ where $k_m$ is the number of  letters of  $E$  
appearing in
$w$ (with or without +) exactly
$m$ times.
   \end{lemma} 
\begin{proof}  
If $(\psi, m)\in \Ker p $, then  $m=0$  and  $  \psi
w=w$.   Since $w(\hat n)=E$, we have   
$\psi=\id$.  Hence $p$ is injective.  

Claim (ii) is proven in the same way as Lemma \ref{ledivisibil} with $M=\{k\in \hat n\,\vert\, \card w^{-1}Ê(w(k))=m\}$. 
\end{proof}

\begin{examp}  For  the  word $W=ABAB  $ in the alphabet $E=\{A,B\}$, the
group 
$\Aut(W)=\ZZ/4\ZZ$  is generated by   $(\psi:E\to E, m=1)$ where $\psi$ permutes $A$ and
$B$.  For   $W= A^+B^+AB $,  
we have $\Aut(W)=1$.      For  the  word
$W=ABACAD$ in the   alphabet $E=\{A,B,C,D\}$, the group 
$\Aut(W)=\ZZ/3\ZZ$ is generated by    $(\psi:E\to E, m=2)$ where $\psi$ fixes $A$ and sends  $B,C,D$ to $C,D,B$,
respectively. 
\end{examp}

\subsection{Words and charts}\label{alphcha}  A chart  $(n, t  )$   gives  rise to 
a   word $W(t)=(w,S)$ as follows. If $n=0$, then $W(t)=\emptyset$.  For $n\geq 1$, the mapping $w$
is the composition of the inclusion $\hat n \hookrightarrow \overline n$ with the   projection 
 $ \overline n\to \overline n /t$.  The set $S\subset \hat n$ consists  of all $k\in \hat n$ such that
$t(k)<0$.   Then $W(t)=(w,S)$ is a  full  word     in the alphabet 
$\overline n /t$.  
It is clear that   $W(\sigma_n t (\sigma_n)^{-1})$  is obtained from  $W(t)$ by the circular permutation.

 \begin{lemma} \label{le7}     For any  chart  $(n, t )$, there is a canonical  group
injection $\Aut (t) \hookrightarrow \Aut(W(t))$.
   \end{lemma} 
\begin{proof} For  $n=0$, both groups   are trivial.  Let $n\geq 1$ and   $ W(t)=(w , S )$. Pick 
$\varphi\in \Aut (t)$.  Recall that  $\varphi=\sigma^m$ where  
$\sigma= \sigma_n:\overline n
\to
\overline n$ and  
$m\in \ZZ/n\ZZ$.   For $k\in \overline  n$, denote its $t$-orbit by $[k]$.  The  equality $\varphi t= t
\varphi$ implies that if
$k, l\in
\overline n$ lie in the same orbit of
$t$, that is if $[k]=[l]$, then
$[\varphi (k)]=[ \varphi (l)]  $. Thus $\varphi$ induces a map   $\psi:  \overline
\pi/t\to \overline
\pi/t$.  Since $\varphi$ is onto, so is $\psi$. Hence $\psi $ is a bijection.  We verify that $\psi w= w\sigma^m$:
 for   $k\in \hat n$,  
$$\psi(w(k))= \psi ([k])= [ \varphi(k)]= [\sigma^m ( k)]=w(\sigma^m ( k)).$$   
The equality $t \sigma^m (k)=\sigma^m t(k)$ implies that $t \sigma^m (k)<0$ iff
$t(k)
<0$. Thus 
   $\sigma^m(S)= S$ and 
$(\psi, m)\in  \Auut (W)$.  The    homomorphism
$\Aut (t)
\to
\Aut (W)$,   $\varphi\mapsto (\psi, m)$    is    injective: if $m=0 $, then 
$\varphi=\sigma^0=\id$. 
\end{proof}

\begin{examp}\label{excicic} The   charts $t_1,t_2,t_3:\overline 4 \to \overline 4$ defined in
Section \ref{automo}    have the same orbits $A=\{-3,-1,1,3\}$ and
$B=\{-4,-2,2,4\}$.  We have   $W(t_1) =A^+B^+ AB$   and $ W(t_2)= W(t_3)=ABAB  $.    Thus,
  different charts  may  yield  the same word   and    in general  
  $\Aut (t)\neq \Aut W(t)$.   \end{examp}

 \subsection{Words associated with curves}\label{char} Any  pointed curve $f $ gives rise to a word $W(f)$ in the alphabet $\cro 
(f)$.
If $f$ is a trivial curve, then $W(f)=\emptyset$.  For a non-trivial $f$, the word $W(f)$ is obtained by first taking the chart  
$(n,t)$  of $f$ and then   taking the associated  word  in the alphabet $\overline n/t=\cro (f)$.    By the results above,  
$$ \Aut  (f) =\Aut (t) \subset \Aut (W (f)) \subset \ZZ/n\ZZ.$$ 
It is easy to read the word $W(f)$ directly from  $f$.  Label  all points of
$\cro (f)$      with  distinct letters (the resulting set of letters is identified with   $ \overline n/t=\cro (f)$).  Label each point
$x\in
\Cross (f) $ with the letter labelling
$f(x)$.  Traverse 
$S^1$   counterclockwise starting from the base point in $S^1-\Cross(f)$ and   write down  consecutively the   letters appearing
when we cross   $\Cross (f)$.   Moreover, crossing a point  $x\in \Cross (f)$  provide the corresponding letter with the superscript $+$
if the flag obtained from  $(x,+)$ by the  flag rotation   is outgoing.
This gives   $W(f) $.  The word   $W(f)$ generalizes the Gauss   word of  a generic curve.    

Similar  constructions apply to   non-pointed curves. Their charts   are  defined    up to conjugation by
circular permutations and    their  words are defined   up to
circular permutations.  

 Example \ref{excicic}
yields different curves with the same associated word.  We now address the realization problem for words.
We say that  a  word $W$  is {\it realized} by a  curve
$f $   if $W$ is congruent to $W(f)$.

	\begin{theor}\label{th222} Every full  word  can be realised by 
 a  pointed   filling 
   curve. 
   \end{theor} 
\begin{proof}  Let 
$W= (w, S )$ be a full  word  of length $n$ in an alphabet $E$.  It suffices to realize $W$ as the word of a chart.  Pick a $w$-equivalence
class
$A\subset
\hat n$.  Let
$a_1,...,a_q$ be  the elements of $A- (A\cap S)$  numerated in an arbitrary way. Let $a_{q+1},...,a_{q+r}$  be  
the elements of $ A\cap S $.
We  cyclically order the  set $\pm A= \{\pm a\}_{ a\in A}  $ by
$$a_1\prec a_2 \prec...   \prec a_{q+1} \prec -a_{q+1} \prec  a_{q+2}\prec -a_{q+2}\prec...\prec  a_{q+r}\prec
- a_{q+r}\prec - a_{q} \prec - a_{q-1}\prec ...   \prec - a_{1}\prec a_1.$$
Let  $t: \pm A\to \pm A$  be the map  sending  each element   to  its immediate follower. 
Applying this procedure to all  $w$-equivalence classes $A\subset 
 \hat n$ we obtain a chart $(n,t: \overline n\to \overline n)$    with orbits  $\{\pm A\}_A$.
  The set $\overline n/t$ can be identified with the set of $w$-equivalence classes in $\hat n$, that is with $E$.
 Under these identifications 
$W(t)=W$. \end{proof}

 \subsection{Words of  self-transversal curves}\label{charseltr} Similarly to charts,  each semichart $(n,v,S)$    gives 
rise to  a   word    $W(v,S) = (w:\hat n\to \hat n/v, S)$ where $w$ is the natural projection from $\hat n$ to the set of $v$-orbits. 
We have  $W(v,S)=W(t)$  where $(n,t)$ is the  straight chart   determined by $(n,v,S)$.
Therefore for a    
self-transversal pointed curve $f $ with chart $(n,t)$ and semichart $(n,v,S)$, we have 
$W(f)  =W(t)=W(v,S)$.    The definition of a semichart implies that the word $W(v,S )=(w,S)$     is {\it odd} in the following sense:  the
intersection of $S$ with any
$w$-equivalence class in
$\hat n$ has an odd number of elements.  

\subsection{Words of coherent curves}\label{wordscharseltr}  We show   that   odd words  bijectively correspond to  coherent
curves.
\begin{theor}\label{th222e} Any  odd full  word     can be realised by 
 a  pointed  coherent   filling 
   curve. This curve is unique up to homeomorphism.   
   \end{theor} 
\begin{proof}  By Theorem \ref{t5pontediksyyy}, it suffices to show that any odd full word $W= (w, S\subset
\hat n)$    arises from a unique coherent semichart.
We define  a permutation $ v=v_w:\hat n\to \hat n$ as follows.  
  For   $k\in \hat n$, let $v (k)$   be the minimal   element of the set $\{ k+1, k+2, ...,n\}$
that is $w$-equivalent to $k$.  If
there are no elements in the latter set $w$-equivalent to $k$,    let  $v (k)$   be  the minimal element of
the set $\{1,2,..., k\}$ that is  $w$-equivalent to $k$. In particular, $v (k)=k$ iff  $w^{-1} (w (k))=\{k\}$.  The resulting mapping $v:
\hat n \to \hat n$  is bijective. Its orbits   are  the  
$w$-equivalence classes in $\hat n$.   The  triple $(n,v ,S)$ is a coherent semichart whose
associated word is $W$. It  follows from the definitions that this is the only such semichart.
 \end{proof}

\begin{corol}\label{t12658745ddddq} The formula $f\mapsto W(f)$ defines a bijective correspondence between  pointed
(resp.\ non-pointed) coherent  filling   curves  considered up to homeomorphism  and  odd full words  
  considered up to congruence (resp.\  considered up to congruence and conjugation by the   circular permutation).
   \end{corol}  

For an odd full word $W$ in an alphabet $E$, the crossings and corners of the corresponding   coherent curve  are labeled (in a
1-to-1 way)  by elements of
$E$. This set of points  has a natural order: the point labeled with $a\in E$ preceeds the point labeled with $b\in E$ if the letter $a$
appears in  
$W$ before
$b$. This order of course is not preserved under circular permutations of $W$.

 The next theorem shows that a coherent filling curve has as many  symmetries as its associated word.

	\begin{theor}\label{th222edfsgsg} For any (non-pointed) coherent   filling 
   curve $f$, we have  $\Aut(f)=\Aut(W(f))$. 
   \end{theor} 
\begin{proof}  Let $(n,v,S)$ be the semichart of $f$ and $W =(w,S)$ be the associated word.  
Since  $(n,v,S)$ is coherent, the map  $v=v_w:\hat n \to \hat n$ is computed from $w$ as in the proof of the previous
theorem. As we know,
$\Aut(f)=\Aut(v,S)\subset \Auut (W)$.  We have to verify  that every   automorphism $(\psi, m\in \ZZ/n\ZZ)$  of $W$ lies in
  $\Aut (v,S)$.  It suffices to show that  $m   $
belongs to the image of $\Aut(v,S)$ under the inclusions 
$\Aut(v,S)\subset
\Aut (W)\subset \ZZ/n\ZZ$.  Thus, we need to check that  the $m$-th power of the circular permutation
$\sigma  :\hat n\to   \hat n$ commutes with   $v$ and  keeps $S$ setwise.
  The last condition follows from the definition of an automorphism of $W$.    The equality   $ 
\psi w=  w \sigma^m $ implies that if    $k,l\in \hat n$ are $w$-equivalent then $\sigma^m(k), \sigma^m(l)$ are
$w$-equivalent.  Thus the map
$\sigma^m:\hat n\to
\hat n$ sends   $w$-equivalence  classes    to $w$-equivalence classes.  Since   this map also preserves the cyclic order in
$\hat n$,  it must  commute  with
$v$.
 \end{proof} 
 
Theorem \ref{th222e}  yields for any  odd word $W=(w , S\subset \hat n)$  a  coherent filling curve    on a
  surface. The genus of this surface   denoted $g(W)$  is a fundamental  geometric  invariant  of  $W$. We
give   an explicit  formula for $g(W)$.

 We shall use the symbol $\prec$ for the standard cyclic order on $\hat n$. Thus for $k, l, m\in \hat n$,  we have $k \prec
l \prec  m$   if   $l\neq k, l\neq m$ and   increasing $k$ by $1,2,...$, we obtain  first $l (\modu n)$ and then $m (\modu n)$.
For $k,l\in \hat n$, set $\vert k, l\vert =\{r\in \hat n\,\vert \, k\prec r \prec l\}$.  In particular, 
$\vert k, k\vert =    \hat n- \{ k\}$.

For  $k\in
\hat n$, set
$d_S(k)=1
$ if
$k\in
\hat n- S$ and
$d_S(k)=0$ if
$k\in S$.  Set    $k^+=v_w (k)$, where $v_w:\hat n \to \hat n$ is the bijection    in the proof of Theorem
\ref{th222e}. 
 For $k,l\in \hat n$, set    $\langle k, l\rangle=1$ if $k\neq l$ and $k$ is $w$-equivalent to $l $; in all the other cases  $\langle k,
l\rangle=0$.   Set 
$$\Delta^0_{k,l}= 
\left\{\begin{array}{ll}
 d_S(k) +d_S(l)+1   , & \mbox{if    $k\prec l \prec k^+ \prec l^+\prec  k$ or $l\prec k \prec l^+ \prec k^+\prec  l$,}\\
\langle k, l\rangle \, d_S(k)\,  d_S (l),   & \mbox{otherwise.}
\end{array} \right.
$$
Finally, define a residue $W_{k,l} \in \ZZ/2\ZZ$   by 
\begin{equation}\label{mainmain} W_{k,l}= \sum_{q\in \vert k, k^+\vert ,  \,r\in \vert l, l^+\vert} \langle  q, r\rangle
+d_S(k) \sum_{r\in \vert l, l^+\vert,   r\neq k^+  } \langle  k, r\rangle
+d_S(l) \sum_{q\in \vert k, k^+\vert ,   q\neq l^+ } \langle  q,  l\rangle+ \Delta^0_{k,l}\, (\modu 2) .\end{equation}

 	\begin{theor}\label{genuszeroor}  For any odd full word $W $ of length $n$, we have
 $g(W)= (1/2)\, {\rank} (W_{k,l})_{k,l\in \hat n}$ where ${\rank}$ is the usual rank of a square $(n\times
n)$-matrix over  the field  $\ZZ/2\ZZ$.
   \end{theor} 

 Theorem \ref{genuszeroor} will be  proven in Section \ref{evopa}.   Note that $W_{k,k}=0$ for all $k\in \hat n$ and $W_{k,l}=d_S(k)\, 
d_S (l)$  for distinct
$w$-equivalent
$k,l\in
\hat n$.  If the alphabet consists of only one letter, all elements of $\hat n$ are $w$-equivalent and 
  $g(W)=(n-\card (S))/2$ for odd $n$ and $g(W)=(n-\card (S)-1)/2$ for even $n\geq 2$.

We call   $W$ {\it planar} if $g(W)=0$, that is if $W$   is realized by a coherent   curve on $S^2$.

 	\begin{corol}\label{genuszeroorert}  An odd full word  $W=(w , S\subset \hat n)$   is planar if and only if 
$W_{k,l}=  0  $ for all $k,l\in \hat n$,
   \end{corol} 

In particular,  if $W$ is planar, then $d_S(k)\,  d_S (l)=0$ for any distinct
$w$-equivalent
$k,l\in
\hat n$. This means that each letter  of the alphabet appears in $W$ with at most one superscipt $+$. As an exercise, the reader may
verify that  $g(A^+A^+A)=1$ and draw   coherent curves on $\RR^2$ representing the words $AAAA^+, ABA^+, BAA^+$.

 \subsection{Example}\label{fchrsitof} Pick relatively prime positive integers $p,n$ with $p\leq n$. The associated {\it Christoffel word} 
in the 1-letter alphabet $\{A\}$   is the word $W=(w:\hat n \to \{A\}, S)$  where $S=\{i\in \hat n\,\vert\, [ip/n]= [(i-1)
p/n ]+1\}$. Clearly, $\card (S)=p$ so that $W$ is odd iff $p$ is odd. If $W$ is odd, then 
$g(W)=(n-p)/2$ for odd $n$ and $g(W)=(n-p-1)/2$ for even $n$.

  \subsection{Remarks}\label{fi:g44refplog}  1.  Deep geometric objects are hiding behind   odd words.   
 As we now know, an   odd full word $W$ in an alphabet $E$ determines  a coherent  (pointed) filling curve $f=f(W)$ on a closed oriented
surface
$\Sigma$ of genus
$g=g(W)$. This     gives rise to an algebraic curve over
$\overline {\QQ}$ and to a    point of the moduli space $\mathcal M_{g,k}$ where
$k=\card (E)$, cf. Introduction and Section
\ref{encoreeee}.    
The   curve $f$ also gives rise to an oriented knot in the total space of the tangent circle bundle  
of
  $\Sigma$.  These constructions do not use the base point of $f$ and therefore the resulting geometric objects are preserved under   
circular permutations of $W$.   

There is a  construction  of a knot from
$f=f(W)$ using its base point. The part of
$f$ lying outside of a small open disc
$D\subset
\Sigma$ surrounding this point   is a proper immersed  interval on $\Sigma-D$. Proceeding as in  \cite{ac}, we can derive from this
immersed interval a knot in a connected sum of $2g$ copies of
$S^1\times S^2$. 

2.  Under the correspondence of Corollary \ref{t12658745ddddq},  the beaming  coherent  filling  curves correspond  
to odd full words $(w,S)$ such that $S$ meets each $w$-equivalence class in   one element. The perfect  
curves  correspond  
to odd full words $(w,S)$ such that $S$ meets each $w$-equivalence class   in its maximal element.

3. A  {\it  Gauss word} is   a  word  $W =(w:\hat n \to E,S)$ such that for all  
$e\in E$,   the set $w^{-1} (e)  $ has  2 elements  and precisely one of them, denoted $e^-$, belongs to $S$.    Such $W$ is
odd and corresponds to  a generic curve 
 with no corners.   The word
$W$  induces a  partition of   $E $ into two disjoint (possibly empty) subsets:  letters 
$e_1,e_2\in E$ belong to the same subset   iff  $e_1^- -e_2^- =0\, (\modu 2)$.
One can show that  whether
$W$ is planar or not depends only on   $w$ and this  partition    of   $E $. Modulo this observation,
Corollary
\ref{genuszeroorert} for   Gauss words is   equivalent to a classical theorem  of Rosenstiehl \cite{ro},
\cite{rr}, see also \cite{ce}.

4. Any  (non-signed) full word 
$W=(w,S =\hat n)$  can be realized by  a     filling  curve $f$ with  $\Auut  (f)=\Aut (W)$.  Indeed, define  a chart $t:\overline n \to
\overline n
$    by 
$t(k)=-k, t(-k)=v_w (k)$ for    $k\in \hat n$. Then $\Auut(t)=\Aut (W)$.

    \section{The genus}\label{evopa}

\subsection{Genus of a chart}\label{greac}  Recall that each closed connected (oriented) surface $\Sigma$ is obtained from the 2-sphere
$S^2$ by attaching several 1-handles.  The number of these
1-handles is   the genus of
$\Sigma$.  
We define the  {\it genus}   $g(t)$  of a chart $(n, t )$  to be  the genus of a    surface  containing a   filling 
 curve with this chart.  The   trivial chart  has genus 
 $0$. 
Suppose   that  $n\geq 1$.
To   compute $g(t)$, recall  the circular  permutation $\sigma=\sigma_n:\overline n \to \overline n$.
Define a map 
$\theta: \overline n \to \overline n$ by $\theta(\pm k)=\mp \sigma^{\pm 1} (k)$ for   $k\in \hat n=\{1,2,...,n\}$.
It is easy to check that $\theta^2=\id$.

 	\begin{theor}\label{genusle}  
 $g(t)=1+ ({n-\card (\overline n/t) - \card (\overline
n/t\theta)})/{2}$.
   \end{theor} 
 \begin{proof}  Let $f:S^1\to \Sigma$ be a filling curve  with chart $(n,t)$.  This curve gives rise to a CW-decomposition
$X$ of $\Sigma$. Its  0-cells are the points of $\cro(f)$, its 1-cells are the components of $f(S^1)-\cro (f)$, and
its  2-cells are the components of $\Sigma - f(S^1)$. Let $\alpha_i$ be the number of $i$-cells of $X$.  Clearly, $\alpha_0=\card
(\overline n/t)$ and 
$\alpha_1=n(f)=n$.  We claim that $\alpha_2=\card (\overline n/t\theta)$. To see this,  recall the flag rotation $t:\Fl(f)\to
\Fl (f)$ and  define  an involution  
$\Theta:\Fl(f)\to
\Fl (f)$
  as follows. Let $x\in \Cross  (f)\subset S^1$ and $\varepsilon =\pm$. Starting at $x$,   move along $S^1$  
counterclockwise if $\varepsilon =+$ and clockwise if $\varepsilon =-$. Let $y$ be the first encountered point of $\Cross
(f)$. Then  $\Theta (x, \varepsilon)= (y, -\varepsilon)$.  Observe that every flag $r\in \Fl(f)$ gives rise to a 2-cell  $D(r)$ of $X$ such 
that  rotating
$r$ around its root in the negative direction until meeting $t^{-1}(r)\in \Fl (f)$ we sweep a subarea of $D(r)$.  It is easy to see   that $D
(t\Theta (r))= D(r)$. Moreover, two flags $r_1,r_2\in \Fl (f)$ verify $D(r_1)=D(r_2)$ if and only if there is
a power of   $t\Theta:\Fl (f) \to \Fl (f)$ transforming $r_1$ into $r_2$.  Thus 
there is a bijective
correspondence between the 2-cells of $X$ and the orbits of   $t\Theta$. 
Under the  identification 
$\Fl (f)=\overline n $ the involution  $\Theta$ is transformed into $\theta$.  Therefore $\alpha_2=\card
(\overline n/t\theta)$. Substituting the expressions for  $\alpha_0, \alpha_1 ,  \alpha_2$ in   the  formula for the Euler
characteristic 
$\alpha_0-\alpha_1 +
\alpha_2 =2-2 g (t)$, we obtain the desired formula for $g(t)$. \end{proof}

\begin{corol} For any chart $(n,t)$, we have $\card (\overline n/t) + \card (\overline
n/t\theta)\leq n+2$. This inequality becomes an  equality   if and only if $(n,t)$ is the chart of a curve on $S^2$. \end{corol}

 Theorem \ref{genusle} allows us to compute     the genus of  specific  charts.   For the   charts
$t_1,t_2,t_3:\overline 4
\to
\overline 4$ from Section \ref{automo},   $g(t_1)=g(t_3)=1$ and $g(t_2)=0$. 
For the  straight chart  $t  :\overline 4 \to \overline 4$
from Section \ref{restc},   $g(t)=1$.

\subsection{Homological formula for the genus}\label{ga20c}  Let  $(n,t )$ be a    chart with
$n\geq 1$. We   compute its genus  $g(t)$  in 
homological terms.  We shall use notation introduced in the proof of Lemma \ref{t2}:  the points 
$x_k=\exp (2k \pi i /n)\in
S^1$ with  $k\in \hat n$,  the graph $\Gamma $, and the
projection $p:S^1\to \Gamma$.  We first specify a set   
 of generators  $ \{h_k\}_{k\in \hat n}$ for the group $H_1(\Gamma)=H_1(\Gamma; \ZZ)$.
 Pick 
$k\in \hat n$. The point $v= p(x_k) $ is a
vertex of $\Gamma$.
Consider the path  
$\gamma_k:[0,1]\to S^1$ starting at
$x_k$ and moving along
$S^1$ counterclockwise until hitting  a point of
$p^{-1} (v)$ for the first time.  (If $p^{-1} (v)=x_k$, then $\gamma_k$   traverses the whole
circle.)  The loop  
$p \gamma_k:[0,1]\to \Gamma$   represents a homology class $h_k\in H_1(\Gamma)$. 
An    induction on the number of vertices of $\Gamma$  shows that 
$H_1(\Gamma)$ is generated by $ \{h_k\}_{k\in \hat n}$.  (It is easy to describe generating relations, but we do not need this.)

Let $\Sigma $ be the closed surface  of genus $g=g(t)$ constructed in  Lemma \ref{t2} by thickening $\Gamma$
  and   gluing 2-disks.  
The inclusion homomorphism  $ H_1(\Gamma)\to 
H_1(\Sigma
)$ is surjective so that  the set $  \{h_k\}_{k\in \hat n}$ generates $H_1(\Sigma)=\ZZ^{2g}$.  The  orientation of $\Sigma$
determines a unimodular   bilinear intersection form 
$B: 
H_1(\Sigma)\times H_1(\Sigma) \to \ZZ$. 
Hence 
		 $$g(t) = (1/2)\, {\rank}  B= (1/2)\, {\rank} (B(h_k, h_l))_{k,l\in \hat n}.$$ In particular, 
  $(n,t)$ is the chart of a curve on $S^2$   if and only if $B(h_k, h_l) =0$ for all $k,l$.

  We give now an explicit formula for the integer  $B(h_k, h_l)$. We  begin with
notation. Recall  the cyclic order induced by $t:\overline n \to \overline n$ on any its orbit (see Section
\ref{prels}). 
For any    $r_1,r_2,r_3 \in \overline n$, we define an integer $\delta (r_1,r_2,r_3 )$: if
  $r_1,r_2,r_3 $ are pairwise distinct and   lie in the same orbit of $t$ 
in  the ($t$-induced) cyclic order $r_1,r_2,r_3, r_1 $, then
$\delta (r_1,r_2,r_3 )=1$; otherwise $\delta (r_1,r_2,r_3 )=0$. Clearly,
$\delta (r_1,r_2,r_3 )= \delta (r_2,r_3,r_1 )=\delta (r_3,r_1,r_2 )$.
For      $r_1,r_2,r_3,r_4\in \overline n$,  set $$\langle  r_1,r_2,r_3,r_4\rangle=
\delta (r_1,r_3,r_2 )\, \delta (r_2,r_4,r_1 )- \delta (r_1,r_4,r_2 )\,  \delta (r_2,r_3,r_1 ) \in \ZZ.$$
Clearly,   $\langle  r_1,r_2,r_3,r_4\rangle=1$ (resp. $\langle  r_1,r_2,r_3,r_4\rangle= -1$)  iff  $r_1,r_2,r_3,r_4$ are pairwise distinct
and   lie in the same $t$-orbit   in  the   cyclic order 
$r_1, r_3, r_2, r_4, r_1$ (resp.   
$r_1, r_4, r_2, r_3, r_1$). 
Otherwise 
$\langle  r_1,r_2,r_3,r_4\rangle=0$. We have  
\begin{equation}\label{anti} \langle  r_1,r_2,r_3,r_4\rangle=-\langle  r_3,r_4,r_1,r_2 \rangle, \end{equation}
$$\langle  r_1,r_2,r_3,r_4\rangle=-\langle  r_2,r_1,r_3,r_4\rangle=-\langle  
r_1,r_2,r_4,r_3\rangle.$$

 For
$k\in
\hat n$, let  $k^+ $ be the minimal   element of the set $\{ k+1, k+2, ...,n\}$ belonging
to the orbit of $t$ containing $k$. 
 If the set $\{ k+1, k+2, ...,n\}$ does not meet the $t$-orbit of   $k$,  then   $k^+$ is   the minimal element of
the set
$\{1,2,..., k\}$ belonging to this orbit.     We shall use the cyclic order $\prec$ on $\hat n$ and the notation   $\vert k, l\vert $
introduced in Section \ref{wordscharseltr}.  For $k,l\in \hat n$, set 
$$\Delta_{k,l}= 
\left\{\begin{array}{ll}
0, & \mbox{if    $k=l$},\\
 \langle -k^+,k,-l^+,l \rangle,  & \mbox{if  $ k\prec k^+\preceq   l\prec l^+\preceq k$}, \\
-\langle  -l^+,l, -l,l^+\rangle, & \mbox{if  $ k\prec l\prec  l^+\prec   k^+$},\\
\langle -k^+, k,-k, k^+ \rangle, & \mbox{if  $ l\prec k \prec k^+\prec l^+ $},\\
\langle -k^+, k,-k, k^+ \rangle-\langle  -l^+,l, -l,l^+\rangle, & \mbox{if      $ k\prec     l^+\preceq  l \prec  k^+\preceq k$},\\ 
 \delta (l,-l, -l^+) -  \delta (k,-k^+, k^+), & \mbox{if  $ k\prec l\prec k^+\prec l^+\prec k $},\\
\delta (l,-l^+, l^+) -\delta (k,-k, -k^+) , & \mbox{if  $l\prec k\prec l^+\prec k^+ \prec l$.}
\end{array} \right.
$$
Note   that    $k,k^+,l,l^+\in \hat n$ satisfy three conditions:  $k=l$ iff
$k^+=l^+$;  if $l=k^+ $, then $k^+\preceq l^+\preceq k$;  if $k=l^+ $, then $l^+\preceq k^+\preceq l$.  The cases listed in
the definition of 
$\Delta_{k,l}$ cover all   possibilities for  such $k,k^+,l,l^+ $.

	\begin{lemma}\label{longl}    For any $k,l \in \hat n$,  
\begin{equation}\label{0001} B(h_k, h_l)=\sum_{q\in \vert k, k^+\vert ,  \,r\in \vert l, l^+\vert} \langle  -q,q, -r,r\rangle
+\sum_{r\in \vert l, l^+\vert} \langle  -k^+,k, -r,r\rangle
+\sum_{q\in \vert k, k^+\vert  } \langle  -q,q, -l^+,l\rangle+ \Delta_{k,l}  .\end{equation} \end{lemma} 

\begin{proof} To compute the   intersection number $B(h,h')$ of
homology classes $h,h'\in H_1(\Sigma)$,  one   presents them  by    transversal  loops $f,f':S^1\to  \Sigma$ such 
that the  set  of common points $f(S^1) \cap f'(S^1)$  consists only of points of multiplicity 1 on both $f$ and $f'$.     Then one assigns
to each    common point  $+1$ if $f$ crosses   $f'$ from left to right at this point and  
$-1$  if $f$ crosses   $f'$ from right to left.  Then  $B(h,h')=f\cdot f'$ is sum of these $\pm 1's$.  

Recall the path $\gamma_k:[0,1]\to S^1$ defined in Section \ref{ga20c}.    The image of
$\gamma_k$ is an arc on
$S^1$ with endpoints
$\gamma_k(0)= x_k$ and $   \gamma_k(1)=x_{k^+} $.  We  denote this arc  by the same symbol $\gamma_k$ and denote the loop
$p\vert_{\gamma_k}: \gamma_k \to \Gamma  \subset \Sigma$ by $f_k$.  By definition, 
$h_k=[f_k]\in H_1(\Sigma)$    where the square brackets denote the homology class of a loop.

  We   prove Formula
\ref{0001} case by case. Throughout the proof we denote  the first, the second, and  the  third  summands on the right-hand side of
(\ref{0001}) by
$(I)_{k,l}, (II)_{k,l}$,   $(III)_{k,l}$  or shorter by $(I), (II)$,  
$ (III)$.  
 
 \vskip0.3truecm (i) If $k=l$, then
$B(h_k,h_l)=\Delta_{k,l}=0$. We must show that  $(I) +(II)+ (III)=0$. By Formula \ref{anti},  $\langle
-q,q,-q,q\rangle =0$   and $ \langle  -q,q, -r,r\rangle
+ \langle  -r,r, -q,q\rangle=0$ for    $q, r\in \hat n$. Hence $(I)= (II)+ (III)=0$.

\vskip0.3truecm (ii) Suppose that  $k\prec k^+\preceq   l\prec l^+\preceq k$. Then   $\gamma_k, \gamma_l  $ do not
meet  except possibly at the endpoints.  The loops $f_k, f_l$ have   a finite set of common points.  We analyse separately 4
possible types of  common points.

  Each pair $q\in \vert k,k^+ \vert ,    r \in \vert l,l^+ \vert $ with $p(x_q)=p(x_r)$ gives   a common point of
$f_k$ and $f_l$. The branch of $f_k$ at $x_q$   is formed by two small arcs representing  the
 incoming  flag $(x_q,-) $ and the   outgoing  flag $(x_q, +) $. Under the identification
$\Fl (f)=\overline n$, these flags  correspond  to $-q$ and $q$, respectively.  Pick a  neighborhood $V$ of
$p(x_q)=p(x_r)$ not containing other vertices of $\Gamma$. If
$\langle  -q,q, -r,r\rangle=0$, then  a  small deformation of
$f_k$ in $V$   makes
$f_k$
 and $f_l$ disjoint  in $V$.   If $\langle  -q,q, -r,r\rangle=\pm 1$, then after a small deformation of $f_k$ in $V$ this loop meets $f_l$
transversally in one point whose  sign  is   $\langle  -q,q, -r,r\rangle$.  Note also that if  $q\in \vert k,k^+ \vert ,    r \in \vert l,l^+ \vert $
and    $p(x_q)\neq p(x_r)$, then $q, r$ do not lie in the same orbit of $t$ and 
$\langle  -q,q, -r,r\rangle=0$.  Thus the    total
contribution to $B(h_k,h_l)$  of the pairs $q\in \vert k,k^+ \vert ,    r \in \vert l,l^+ \vert $ with $p(x_q)=p(x_r)$  is equal
to $(I)$.

  Each  $ r \in \vert l,l^+ \vert  $ with $p(x_k)=p(x_r)$ gives   a common point of
$f_k$ and $f_l$. The branch of $f_k$ at  $ x_k$ is formed by two   arcs representing  the
   flags $(x_{k^+},-) $ and  $(x_k, +) $. Under the
identification
$\Fl (f)=\overline n$, they  correspond  to $-k^+$ and $k$.  As above,  the
contribution of  this   point to $B(h_k,h_l)$  is   $\langle  -k^+,k, -r,r\rangle$.      Similarly, each  $q\in \vert k,k^+ \vert  $ with
$p(x_q)=p(x_l)$ contributes $ \langle  -q,q, -l^+,l\rangle$. This gives  (II) and    (III). Finally,  in the case where
$p(x_k)=p(x_l)$,   this point contributes $ \langle -k^+,k,-l^+,l \rangle=\Delta_{k,l}$.

\vskip0.3truecm (iii) Suppose that  $k\prec l\prec  l^+\prec   k^+$.  The  points $x_l,
x_{l^+}$ split  
$\gamma_k$ into three subarcs: an arc $\alpha$ leading from $x_k$ to
$x_l$, the arc $\gamma_l$ leading from $x_l$ to
$x_{l^+}$ and an arc $\beta$  leading from  
$x_{l^+}$    to
$x_{k^+}$.   Restricting 
$p$ to   $\alpha$ and $\beta$ we obtain two composable paths in $\Sigma$ whose composition is a loop,   $f$.
It is clear that $h_k=h_l+ [f]$. Since $B$ is skew-symmetric,  
$B(h_k,h_l)=B([f], h_l)$.  Since  the loops $f$ and $f_l$ meet only at a finite set 
of points, the same argument as in   (ii) computes $B(h_k,h_l)=B([f], h_l)$ to be  
$$ \sum_{q\in \vert k, l \vert  \cup \vert l^+, k^+\vert, \,  r \in \vert l,l^+ \vert  }
\langle  -q,q, -r,r\rangle + \sum_{r\in \vert l, l^+\vert} \langle -k^+, k, -r,
r\rangle
 +
\sum_{q\in  \vert k, l \vert  \cup \vert l^+, k^+\vert } \langle  -q,q, -l^+,l\rangle  + \langle  -l, l^+, -l^+,l\rangle
.$$ Here   $p(x_k)\neq p(x_l)$ and $p(x_r)\neq p(x_l)$ for all $r\in \vert l, l^+\vert$ because the set $\vert k,
k^+
\vert$ does not  meet  the $t$-orbit of $k$ and   the set $\vert l,
l^+
\vert$ does not  meet  the $t$-orbit of $l$. Denote the   three big sums in this expression for  $B(h_k,h_l)$  by $(I)', (II)', (III)'$,
respectively.  It is clear that $(II)'=(II)$ and
$$(I)=(I)'+ \sum_{q\in \vert  l , l^+ \vert  , \,  r \in \vert l,l^+ \vert  }
\langle  -q,q, -r,r\rangle + \sum_{q\in  \{l,l^+\} , \,  r \in \vert l,l^+ \vert  }
\langle  -q,q, -r,r\rangle= (I)'+0+0=(I)'.$$
  Similarly,    $(III)'=(III)$.  It
remains to observe that $\langle  -l, l^+, -l^+,l\rangle=-\langle  -l^+,l, -l, l^+\rangle=\Delta_{k,l}$.

\vskip0.3truecm (iv) If  $ l\prec k \prec k^+\prec l^+$, then by (iii), 
$$B(h_k, h_l)=-B(h_l, h_k) =-(I)_{l,k} -(II)_{l,k} -(III)_{l,k}- \Delta_{l,k}.$$
Formula \ref{anti} implies that  $-(I)_{l,k}=(I)_{k,l}$, $-(II)_{l,k}=(III)_{k,l}$, $-(III)_{l,k}=(II)_{k,l}$, and
$ - \Delta_{l,k}=  \Delta_{k,l}$.

\vskip0.3truecm (v).   Suppose that  $ k\prec     l^+\preceq    l \prec  k ^+\preceq  k$.  If $k=k^+\neq  l=l^+$, then   $h_k=h_l= [p]$
and
$B(h_k,h_l)=0$. The right hand side of  (\ref{0001}) is  
easily computed to be 
$\sum_{q,   r\in \hat n} \langle  -q,q, -r,r\rangle=0$ (one should use that $ \langle -k,k, -l, l \rangle=0$ which is due to the
fact that   $k=k^+$   is the only element in its $t$-orbit).

Assume that  $ k\prec     l^+\prec    l \prec  k ^+\prec  k$.   (We leave the    cases
$ k\prec     l^+=  l \prec  k ^+\prec k$ and $ k\prec     l^+\prec   l \prec  k ^+= k$ to the reader.)     
The circle 
$S^1$ splits  into four arcs
$\alpha,
\beta,
\alpha', \gamma $ leading from $x_k$ to $x_{l^+}$, from
$x_{l^+}$ to $x_l$, from $x_l$ to $x_{k^+}$, and from $x_{k^+}$ to $x_k$, respectively. 
Restricting 
$p:S^1\to \Gamma\subset \Sigma$ to   $\alpha$ and $\alpha'$ we obtain two composable paths   whose composition is a
loop,   
$f$. Restricting 
$p$ to   $\beta$ (resp. $\gamma$) we obtain a loop $g$ (resp.   $h$)  in $\Sigma$. Clearly, $h_k=[f]+[g]$ and $h_l=[f]+[h]$.
Therefore
$B(h_k,h_l)= B(h_k, [h]) + B ([g], [f])$. The arcs   
$\gamma_k= \alpha \cup \beta \cup \alpha'$ and $\gamma $ have the same endpoints and are otherwise disjoint. 
The same argument as in   (ii) gives 
\begin{equation}\label{000777} B(h_k, [h])=\sum_{q\in \vert k, k^+\vert , \,  r \in \vert k^+,k \vert  }
\langle  -q,q, -r,r\rangle+\sum_{r \in \vert k^+,k \vert }
\langle -k^+, k, -r,r\rangle + \langle -k^+, k, -k, k^+\rangle
\end{equation}
where we use    that the set $\vert k, k^+
\vert$ does not  meet  the $t$-orbit of $k$.   The  same method   gives  
\begin{equation}\label{00012} B ([g], [f])= 
\sum_{q\in \vert  l^+,l \vert,   \,  r \in \vert k,l^+  \vert \cup \vert l, k^+\vert}
\langle  -q,q, -r,r\rangle +\sum_{q \in \vert  l^+,l  \vert} \langle    -q,q, -l^+,l\rangle + \langle      -l, l^+, -l^+,l \rangle.
\end{equation}
 The absence of further summands    is due to the fact that the set $ \vert   l^+, l \vert 
\subset \vert k, k^+\vert   $ does not  meet  the $t$-orbit  of $k$ and  the set  $ \vert k, l^+ \vert 
\cup \vert l, k^+\vert \subset  \vert l, l^+\vert $ does not  meet  the $t$-orbit of  $l$. 
Denote the first and the second summands on the right hand side of (\ref{000777}) (resp. (\ref{00012})) by $(I)'$ and $(II)'$
(resp.   $(I)''$ and $(II)''$). 
Observe   that   $(II)'=  (II) $
since the set $\vert l, l^+  \vert- \vert k^+,k \vert  $    meets  the $t$-orbit of $k$  in two points $r=k,k^+$ and for both    $\langle
-k^+, k, -r,r\rangle=0$. Also   $(II)''=(III)$ since   the set $\vert k, k^+  \vert- \vert l^+,l \vert 
$    meets  the $t$-orbit of $l$   in two points $q=l,l^+$ and for both    $\langle    -q,q, -l^+,l\rangle=0$.
These computations yield  $B(h_k,h_l)  =(I)'+(II)'+(II)+(III)+\Delta_{k,l}$.  We claim that
$ (I)'+(I)''=(I)$.  The equality 
$ \vert l,l^+  \vert=\vert l,k^+  \vert
\cup \vert k^+,k  \vert\cup \vert k,l^+  \vert \cup \{k,k^+\}$ implies  that
$$(I)= (I)'+ \sum_{q\in \vert k, k^+\vert , \,  r \in \vert k,l^+  \vert \cup \vert l, k^+  \vert  }
\langle  -q,q, -r,r\rangle.$$
Using that $ \vert k,k^+  \vert=\vert k,l^+  \vert
\cup \vert l^+,l  \vert\cup \vert l,k^+  \vert \cup \{l, l^+\}$, we    further  expand 
$$(I)= (I)'+  (I)''+ \sum_{q, r\in   \vert k,l^+  \vert \cup \vert l, k^+  \vert  }
\langle  -q,q, -r,r\rangle + \sum_{q\in \{l, l^+\}, r\in   \vert k,l^+  \vert \cup \vert l, k^+  \vert  }
\langle  -q,q, -r,r\rangle=(I)'+  (I)''$$
since both  sums in the  middle term are $0$.

\vskip0.3truecm (vi).  Suppose that $k\prec l\prec k^+\prec l^+\prec k$. Then $p(x_k)\neq p(x_l)$ and both  loops 
$f_k, f_l$ representing $h_k, h_l$   contain the image of the arc on $S^1$ leading from $x_l$ to
$x_{k^+}$.  For $j\in \overline n$, denote a small arc on $\Gamma$ representing the  corresponding flag  by $\alpha_j$.  Pushing 
$f_k$  slightly  to its right in $\Sigma$, we obtain a loop, $f_k^+$,  transversal to $f_l$. We do it so that the 
point 
$p(x_k) $ is pushed  to a point in $\Sigma-\Gamma$ lying   between    $\alpha_k$ and $\alpha_{t^{-1} (k)}$.  Then  
  $B(h_k, h_l) =f^+_k\cdot f_l$.  In a small neighborhood $V$ of   $p(x_l)$ the branch of  $f_l$ at $x_l$ is formed by the incoming arc 
$\alpha_{-l^+} $ and the outgoing arc $\alpha_l $. The loop $f_k$  goes through $p(x_l)$ as many times as there are  elements   $q\in
\vert k,  k^+\vert
$ lying in the $t$-orbit of $l$. (All such $q$ lie in $\vert k, l\vert \cup \{l\}$.) For   $q\neq l$, the corresponding contribution to
$f^+_k\cdot f_l$   is  
$\langle -q, q, -l^+, l \rangle$.  If       $\alpha_{-l^+} $ lies on the left of     $\alpha_{-l} \cup \alpha_{l} $, then the branch of $f_k^+$ at
$x_l$ does not meet the branch of $f_l$ at  $x_l$.  If         $\alpha_{-l^+} $ lies  on the right  of
   $\alpha_{-l} \cup \alpha_{l} $, then the branch of $f_k^+$ at
$x_l$    meets the branch of $f_l$ at  $x_l$  transversally in 1 point, whose sign is   $+1$.
This gives a total of  
$$\sum_{q\in \vert k,  k^+\vert-\{l\} }\langle -q, q, -l^+, l \rangle +\delta (l,-l, -l^+)=\sum_{q\in \vert k,  k^+\vert  }\langle -q, q,
-l^+, l \rangle +\delta (l,-l, -l^+)=(III)+\delta (l,-l, -l^+).$$ 
Similarly, the common points of $f_k^+, f_l$ lying  a small neighborhood  $V'$  of  
$p(x_k)$     contribute 
$$\sum_{  r \in \vert l,l^+ \vert  } \langle  -k^+,k, -r,r\rangle -  \delta (k,-k^+, k^+)=(II)-    \delta (k,-k^+, k^+) .$$
Finally, the   contribution to $f^+_k\cdot f_l$  of the common points of $f_k^+$ and $f_l$ lying  outside of $V\cup V'$  is  
$(I)$.

 \vskip0.3truecm (vii).  The case    $ l\prec k\prec l^+\prec k^+\prec l$ is deduced from the previous
one by  the same argument as in (iv). \end{proof}

\subsection{Proof of Theorem \ref{genuszeroor}}    Since the   form $B$ in Section \ref {ga20c} is unimodular, its rank remains the same after tensor multiplication
by
$\ZZ/2\ZZ$. Therefore to prove Theorem \ref{genuszeroor} it suffices to show that  $W_{k,l}=B(h_k, h_l) \, (\modu 2)$ for all $k,l\in
\hat n$ provided  $(n,t)$ is the chart of the   coherent curve determined by   $W$.  This is a consequence of   the following
equalities modulo 2. For  all 
$k,l,  q,r\in
\hat n$, we have 
$ \langle  -q,q, -r,r\rangle= \langle  q, r\rangle  $. We have  $ \langle  -k^+,k, -r,r\rangle=d_S(k)  \langle   k, r\rangle\,
 $ if
$r\neq k^+$ and 
$ \langle  -k^+,k, -r,r\rangle= 0 $ if $r= k^+$.
Similarly,  $ \langle  -q,q, -l^+,l\rangle= d_S(l) \langle   q,  l\rangle $ if $q\neq l^+$ and 
 $ \langle  -q,q, -l^+,l\rangle= 0 $ if $q= l^+$. Finally, 
$\Delta_{k,l}=\Delta^0_{k,l} \, (\modu 2)$ because    $\langle -k^+,k,-l^+,l \rangle=\langle k, l\rangle \, d_S(k) d_S (l) $,
$\langle  -k^+,k, -k,k^+\rangle= 0$, 
$\delta(k,-k,-k^+)=d_S(k)  $ and 
$\delta(k,-k^+,k^+)=1+d_S(k)  $.

                      \end{document}